\documentclass[12pt]{amsart}
\usepackage{amsmath,amsthm,amsfonts,amscd,amssymb}
\usepackage[all]{xy}

\theoremstyle{plain}

\newtheorem{Thm}{Theorem}[section]
\newtheorem{Lem}[Thm]{Lemma}
\newtheorem{Cor}[Thm]{Corollary}
\newtheorem{Prop}[Thm]{Proposition}
\newtheorem{Def/Prop}[Thm]{Definition/Proposition}
\newtheorem{Rem}[Thm]{Remark}

\newtheorem*{Thm*}{Theorem}
\newtheorem*{Rem*}{Remark}
\newtheorem*{Lem*}{Lemma}

\newtheorem*{Cor*}{Corollary}
\newtheorem*{Prop*}{Proposition}
\theoremstyle{definition}
\newtheorem{Def}[Thm]{Definition}

\newcommand{\ra}{\rightarrow}
\newcommand{\hra}{\hookrightarrow}
\newcommand{\dra}{\dashrightarrow}

\newcommand{\ov}{\overline}
\newcommand{\til}{\widetilde}

\newcommand{\bZ}{\mathbb{Z}}
\newcommand{\bP}{\mathbb{P}}
\newcommand{\bA}{\mathbb{A}}

\newcommand{\cH}{\mathcal{H}}
\newcommand{\cU}{\mathcal{U}}
\newcommand{\cF}{\mathcal{F}}
\newcommand{\cO}{\mathcal{O}}
\newcommand{\cJ}{\mathcal{J}}

\newcommand{\cK}{\mathcal{K}}

\newcommand{\et}{\acute{e}t}

\begin{document}

\title{Zero cycles on homogeneous varieties}

\thanks{This material is based upon work supported by the National
Science Foundation under agreement No. DMS-0111298. Any opinions,
findings and conclusions or recomendations expressed in this material
are those of the author and do not necessarily reflect the views of
the National Science Foundation.}

\author{Daniel Krashen}

\begin{abstract}

In this paper we study the group $A_0(X)$ of zero dimensional
cycles of degree $0$ modulo rational equivalence on a projective
homogeneous algebraic variety $X$. To do this we translate rational
equivalence of $0$-cycles on a projective variety into R-equivalence
on symmetric powers of the variety. For certain homogeneous varieties,
we then relate these symmetric powers to moduli spaces of \'etale
subalgebras of central simple algebras which we construct. This allows
us to show $A_0(X) = 0$ for certain classes of homogeneous
varieties, extending previous results of Swan / Karpenko, of
Merkurjev, and of Panin.

\end{abstract}

\maketitle

\section{Introduction}

The study of algebraic cycles on quadric hypersurfaces has turned out
to be unreasonably successful in its applications to quadratic
forms. Karpenko, Izhboldhin, Rost, Merkurjev, Vishik and Voevodsky, to
name a few, have used and developed the theory of algebraic cycles in
order to solve a number of outstanding conjectures, most notably
Voevodsky's recent proof of the Milnor conjecture.

In part inspired by these great successes, there is much interest in studying
algebraic cycles on and motives of general projective homogeneous
varieties, beyond the quadric hypersurfaces which arise in
applications to quadratic forms. Significant progress has been made by
various authors in this direction (\cite{Kar:Cell, Bro:BB, CGM,
SeZa}).

Despite the progress in understanding general projective homogeneous
varieties, the Chow groups of $0$-dimensional cycles for such
varieties have remained somewhat mysterious. Whereas computations have
been performed in various cases (see for example Swan \cite{Swan} and
Merkurjev \cite{Merk:0cycles}), the topic has so far resisted general
statements or conjectures.

In this paper, we compute the Chow group of zero cycles on various
projective homogeneous varieties by showing that the group $A_0(X)$ of
$0$ dimensional cycles of degree $0$ modulo rational equivalence is
trivial in many cases. We give examples of this for certain
homogeneous varieties for groups of each of the classical types $A_n,
B_n, C_n, D_n$.

More precisely, in the $A_n$ case (theorem \ref{0_gsbv}), we show that
$A_0(X) = 0$ for $X$ a Severi-Brauer variety (recovering a result of
Panin), and for certain cases when $X$ is a Severi-Brauer flag
variety. In all of these examples, we assume that either $F$ is
perfect or $char(F)$ doesn't divide the index of the underlying
central simple algebra.

In the $B_n$ and $D_n$ cases (theorem \ref{0_orth}), we show that
$A_0(X) = 0$ for any (orthogonal) involution variety $X$, assuming
that $char(F) \neq 2$. Involution varieties are twisted forms of
quadric hypersurfaces introduced in \cite{Tao}, and are defined in
section \ref{inv}. This generalizes previous results of Swan
(\cite{Swan}) and Karpenko who proved this when $X$ is a quadric
hypersurface, and Merkurjev (\cite{Merk:0cycles}) who proved this when
$X$ has index $2$ (see section \ref{notation} for the definition of
index).

In the $C_n$ case (theorem \ref{0_symp}), we show that $A_0(X) = 0$
for $X = V_2(A, \sigma)$ a 2'nd generalized involution variety for a
central simple algebra $A$ with symplectic involution $\sigma$ (see
section \ref{inv}) when $ind(X) = 1 or 2$ and $char(F) \neq 2$. This
gives the first nontrivial computations of this group for such
varieties. The case of higher index is still open.

To obtain our results we relate the Chow group of $0$-dimensional
cycles to the more geometrically naive notion of R-equivalence
(i.e. connecting points with rational curves) on symmetric powers of
the original variety, along with the slightly weaker notion of
H-equivalence which we introduce. This is explained in section
\ref{symm}. Although in some sense, this idea is not new - various
aspects of this idea over the complex field appear in \cite{Sam}, and
similar ideas were used in Swan's paper (\cite{Swan}), our formulation
of this principle allows us to more fully exploit its uses.

From here, we show that the symmetric powers of certain homogeneous
varieties may be related to spaces which parametrize commutative
\'etale subalgebras in a central simple algebra. To make this
connection precise, we define moduli spaces of \'etale subalgebras in
section \ref{moduli}. These spaces are very interesting in their own
right, as many open questions in the area of central simple algebras
concern the existence and structure of certain types of subfields in a
division algebra. In sections \ref{maxl}, \ref{deg4} and
\ref{exp2_subsec} we determine show that in certain cases these moduli
spaces are R-trivial, and in sections \ref{gsbv} and \ref{inv} we
apply this to determining the Chow group of zero cycles for certain
homogeneous varieties.

There are various known results concerning the group $A_0(X)$ for
geometrically rationally connected varieties over certain fields,
particularly the finite, local, and global cases (see
\cite{KolSza:RCFF}, \cite{Kol:RCLF}, and \cite{CT}).  For example,
Colliot-Th\'el\`ene has conjectured that the torsion in $CH_0(X)$ is
finitely generated when $F$ is $p$-adic, and has obtained positive
results in certain cases (\cite{CT}).

Over an arbitrary ground field, it is clear that the geometrically
rationally connected varieties may have very complicated groups of
zero cycles, and so it appears difficult to know which classes of
varieties have $A_0(X) = 0$. Even restricting to projective
homogeneous varieties is not sufficient for this. For example,
A. Vishik has pointed out the following example using a result of
Karpenko and Merkurjev (\cite{KarpMerk}):

\begin{Prop}
One may find a field $F$ and a quadratic form $q$ over a vector space
$V/F$ such that if we let $X$ be the variety of $2$-dimensional
totally isotropic subspaces of $V$, the group of $CH_0(X)$ is
infinitely generated (and therefore so is $A_0(X)$).
\end{Prop}
\begin{proof}
For a given quadratic form $q$ on $V/F$ we may construct the variety
$X$ as above. Let $Q$ be the quadric hypersurface in $\mathbb{P}(V)$
defined by the vanishing of $q$. Thinking of points in $Q$ as
isotropic lines in $V$, we may construct a Chow correspondence from
$X$ to $Q$ by setting $Z \in X \times Q$ to be the subvariety
described as
$$\{(x, q) \in X \times Q | q \subset x\}.$$ This defines a
homomorphism $CH_0(X) \ra CH_1(Q)$.

In \cite{KarpMerk}, the authors exhibit a quadratic form $q$ on a $7$
dimensional vector space such that the associated $5$-dimensional
quadric $Q$ has an infinite family of independent nontrivial torsion
cycles $z_i \in A^4(Q) = CH_1(Q)$. One may check by inspection that
these cycles are in the image of the Chow correspondence above, and
therefore give infinitely many independent nontrivial elements in
$CH_0(X)$.
\end{proof}

I am grateful to A. Merkurjev who suggested this problem to me while I
was a VIGRE assistant professor at UCLA, and whose helpful comments on
various drafts of this paper were extremely useful. I would also like
to thank D. Saltman who suggested to me the idea of using Pfaffians to
prove theorem \ref{exp2}, and I. Panin who explained to me how to
concretely think of the varieties associated to symplectic
involutions. I am also grateful for the comments of an anonymous
referee who reccomended the use of Hilbert schemes after reading a
previous version of this paper. The use of Hilbert schemes of points
has considerably cleaned up and shortened the exposition of the paper,
as well as done away with almost all assumptions about the
characteristic of the ground field.

After the appearance of this paper in preprint form, Viktor Petrov,
Nikita Semenov and Kirill Zainoulline have subsequently applied these
methods to compute groups of $0$ cycles on homogeneous varieties for
various exceptional groups (\cite{PSZ}).

\section{Preliminaries and notation} \label{notation}

Let $F$ be a field. All schemes will be assumed to be separated and of
finite type over a field (generally $F$ unless specified
otherwise). By a variety, we mean an integral scheme. If $Z$ is a
closed subscheme of a scheme $X$, we let $[Z]$ denote the
corresponding cycle. Suppose $X$ and $Y$ are schemes over $F$. For an
extension field $L/F$ we denote by $X_L$ the fiber product $X
\times_{Spec(F)} Spec(L)$. For a morphism $f: X \ra Y$, we write $f(L)
: X(L) \ra Y(L)$ for the induced map on the $L$-points. We denote by
$F(X)$ the function field of $X$. We define the index of a scheme $X$,
as
$$ind(X) = GCD\{[L : F] \ | \ L/F \text{ finite field extension and }
X(L) \neq \emptyset\}.$$

If $A$ is a central simple $F$ algebra, we recall that its dimension
is a square, and we define the degree of $A$, $deg(A) =
\sqrt{dim_F(A)}$. We may write such an $A = M_m(D)$ for some division
algebra $D$ unique up to isomorphism, and we define the index of $A$,
$ind(A) = deg(D)$. We let $exp(A)$ denote the order of the class of
$A$ in the Brauer group $Br(F)$. If $M$ is a finite $A$ module, we
follow \cite{BofInv} and define the reduced dimension of $M$ to be
$rdim(M) = dim_F(M) / deg(A)$.

We will make frequent use of symmetric powers and Hilbert schemes of
points. For this purpose, we will make the following notational
shorthands. For a quasiprojective variety $X$ over $F$, we define the
symmetric power $S^n X$ to be the quotient $X^n/ S_n$. We define
$X^n_\circ$ to be the configuration space of $n$ distinct points on
$X$ - i.e. $X^n_\circ = X^n \setminus \Delta$, where $\Delta$ is the
big diagonal. We let $X^{(n)}$ be the quotient $X^n_\circ/S_n$. Note
that the quotient morphism $X^n_\circ \ra X^{(n)}$ is \'etale. For $X$
quasiprojective, we let $X^{[n]}$ denote the Hilbert scheme of $n$
points on $X$ and $\cU^n_X \subset X^{[n]} \times X$ denote the
universal family over the Hilbert scheme $X^{[n]}$. Note that
$X^{(n)}$ is a dense open subscheme of $X^{[n]}$ if $dim(X) \geq 1$.

In the case that $X$ is given as a subscheme of a Grassmannian $X
\subset Gr(k, m)$, we let $X^n_* \subset X^n$ denote the open
subscheme consisting of collections of $n$ subspaces $W_1, \ldots,
W_n$ which are linearly independent, and we set $X^{(n)}_* = X^n_* /
S_n$. 

For a scheme $X$, we define $Z(X)$ to be the set of $0$ dimensional
cycles on $X$ and $Z^n_{eff}(X)$ to the the subset of degree $n$
effective cycles in $Z(X)$. We have a set map $X^{[n]}(F) \ra
Z^n_{eff}(X)$ defined by taking a subscheme $z \subset X$ of degree
$n$ to its fundamental class $[z]$. This gives a bijection between the
cycles which are a disjoint union of spectrums of separable field
extensions of $F$, and points in $X^{(n)}(F) \subset X^{[n]}(F)$. We
will occasionally have to make use of cycles of other dimensions, and
we will use the notation $C_i(X)$ to represent the group of
$i$-dimensional cycles on $X$.

We say that a field $L$ is prime to $p$ closed if every finite
algebraic extension $E/L$ has degree a power of $p$. An algebraic
extension $L/F$ is called a prime to $p$ closure if for every finite
subextension $F \subset L_0 \subset L$, $[L:F]$ is prime to $p$, and
$L$ is prime to $p$-closed.

\begin{Lem} \label{index_hilb_open}
Suppose $X$ is a scheme over $F$ with $ind(X) = n$, where either
$char(F)$ doesn't divide $n$ or $F$ is perfect. Then $X^{(n)}(F) =
X^{[n]}(F)$.
\end{Lem}
\begin{proof}
Given a point $x \in X^{[n]}(F)$, $x$ corresponds to a finite
subscheme $Spec(R) \subset X$, where $R$ is a commutative $F$-algebra
of dimension $n$. By taking a quotient by a maximal ideal of $R$, we
obtain subscheme $Spec(L) \subset X$, $L$ a field of degree at most
$n$. Since $ind(X) = n$, we immediately conclude $Spec(R) = Spec(L)$
and so $R$ is a field. By our hypothesis, $R$ is a separable field
extension, and so we see that $x$ corresponds to a point in
$X^{(n)}(F)$ as claimed.
\end{proof}

\begin{Lem} \label{prime_reduction_new}
Let $X$ be a proper variety such that for any extension field $L/F$,
$X(L) \neq \emptyset$ implies $A_0(X_L) = 0$. If $A_0(X_{F_p}) = 0$
for each prime $p$ dividing $ind(X)$ and every prime to $p$ closure
$F_p/F$ then $A_0(X) = 0$.
\end{Lem}
\begin{proof}
Suppose first that $p$ does not divide $ind(X)$. It then follows that
$X(F_p) \neq \emptyset$, and hence by the hypotheses, $A_0(X_{F_p}) =
0$. Therefore, the conditions of the lemma imply $A_0(X_{F_p}) = 0$
for all $p$.

We will show that $A_0(X) = 0$ by showing that the degree map $deg :
CH_0(X) \ra \bZ$ is injective. Let $deg_p$ be the degree map after
fibering with $F_p$. Consider the natural map $\pi_p: X_{F_p} \ra X$,
which is a flat morphism.  Let $\alpha \in ker(deg)$, and assume that
$deg_p$ is injective. In this case, $\pi_p^*(\alpha) \in ker(deg_p) =
0$. This means that we may find irreducible curves $Z_i \subset
X_{F_p}$, and a rational function $r_i \in R(Z_i)$, such that $\sum
div \ r_i = \alpha$. But since these subvarieties $Z_i$, and functions
$r_i$ involve only a finite number of coefficients, they are defined
over an finite degree intermediate field $E$, where $F \subset E
\subset F_p$. But now we have that if $\pi_E: X_E \ra X$ is the
natural map, then $\pi_E^*\alpha = 0$. But ${\pi_E}_* \pi_E^* \alpha =
[E:F] \alpha$ tells us that $[E:F] \alpha = 0$, and so $[E:F] \in
ann_\bZ(\alpha)$. Therefore, since $[E:F]$ is prime to $p$,
$ann_\bZ(\alpha) \not\in [E:F]$. But because this holds for every
prime $p$, we must have that $ann_\bZ(\alpha)$ is not contained in any
maximal ideal of $\bZ$ and hence $ann_\bZ(\alpha) = \bZ$. But this
implies that $\alpha = 0$.
\end{proof}

\section{Cycles and equivalence relations} \label{symm}

Let $X$ be a scheme. We say that two points $p_1, p_2 \in X(F)$ are
elementarily linked if there exists a rational morphism $\phi :
\mathbb{P}^1 \dra X$ such that $p_1, p_2 \in im(\phi(F))$. We define
R-equivalence to be the equivalence relation generated by this
relation. Let $X(F)/R$ denote the set of equivalence classes of points
in $X(F)$ under $R$-equivalence. We say that $X$ is R-trivial in case
$X(F)/R$ is a set of cardinality $1$. 

If $f : X \ra Y$ is a morphism, we obtain a map of sets $X(F)/R \ra
Y(F)/R$ which we denote by $f_R$. Note that this is well defined,
since if $p, q\in X(F)$ are elementarily linked via a rational map
$\mathbb{P}^1 \dra X$, then the composition $\mathbb{P}^1 \dra X \ra
Y$ shows that $f(p)$ and $f(q)$ are elementarily linked as well. 

Given points $x, y \in X^{[n]}(F)$, we say that $x$ and $y$ are
elementarily H-linked if there is a morphism $\phi : \bP^1 \ra
X^{[n]}$ such that $[\phi(0)] = [x], [\phi(1)] = [y]$. We define
H-equivalence, denoted $\sim_H$, to be the equivalence relation
generated by elementary H-linkage. We say that an open subscheme $U
\subset X^{[n]}$ is H-trivial if the H-equivalence classes $U(F)/H$
form a set with one element. Note that for $x, y \in X^{(n)}$, $[x] =
[y]$ if and only if $x = y$.

We remark that the notions of R and H equivalence carry over in
relative versions for any base scheme $S$ by replacing $\bP^1_F$ with
$\bP^1_S$. In particular, if $S \cong Spec(\oplus E_i)$ where each $E_i$
is a field, it is easy to check that two points are R or H equivalent
if and only if the corresponding points are equivalent with respect to
each $E_i$.

The first lemma we prove gives some justification for considering
H-equivalence:
\begin{Lem} \label{symmch_new}
Suppose $X$ is a projective variety, and $\alpha, \beta \in
X^{[n]}$. If $\alpha$ and $\beta$ are H-equivalent, then $[\alpha]$
and $[\beta]$ are rationally equivalent.
\end{Lem}
\begin{proof}
Without loss of generality, we may assume that $\alpha$ and $\beta$
are elementarily linked, and choose a morphism $\phi : \bP^1 \ra
X^{[n]}$ connecting these points (we may assume $\phi$ is a morphism
and not just a rational map since the Hilbert scheme is
proper). Pulling back the universal family on $X^{[n]}$ along $\phi$,
we obtain a flat family $F \subset X \times \bP^1$ of $0$ dimensional
subvarieties of $X$ of degree $n$ on $\bP^1$. By \cite{Ful:IT},
section 1.6, any two specializations of this to points in $\bP^1$
are rationally equivalent. In particular, $\alpha$ and $\beta$ are
rationally equivalent. 
\end{proof}

\begin{Lem} \label{hilb_chow}
Suppose $X$ is a projective variety over $F$ with $dim(X) \geq
1$. Then the map $X^{[n]}(F) \ra Z^n_{eff}(X)$ is surjective.
\end{Lem}
\begin{proof}
Let $z \subset X$ be an irreducible effective $0$ cycle, say $z \cong
Spec(L)$ for $L/F$ a finite field extension. It suffices to show that
for any $r>1$, there is a subscheme $\til{z} \subset X$ with
$[\til{z}] = r[z]$. Without loss of generality, we may assume that $X
= Spec(R)$ is affine. Let $m \subset R$ be the maximal ideal
corresponding to $z$. Since $dim(R) \geq 1$, we know that
$length(R/m^k)$ is unbounded as $k$ increases. In particular, there
exists $k>0$ such that $R/m^k$ has length $\geq r$ and $R/m^{k-1}$ has
length $< r$. Now consider the module $m^{k-1} / m^k$. We need only
show that this module has a submodule $M$ with $length(M) = r -
length(R/m^{k-1})$. But submodules of $m^{k-1}/m^k$ are the same as
$L$ vector spaces with length corresponding to dimension. Since we
have subspaces of any desired size, we are done.
\end{proof}

With this in mind, it is reasonable to extend the definition of
elementary H-linkage and H-equivalence to cycles. Namely, if $x,y$ are
effective zero cycles of degree $n$ on a regular variety $X$, we say
that they are elementarily H-linked (H-equivalent resp.), if there
exist $x', y' \in X^{[n]}(F)$ with $[x'] = x, [y'] = y$ such that $x'$
and $x'$ are elementarily H-linked (H-equivalent resp.).

\begin{Cor} \label{cycle_hilb_H}
There is a natural bijection $X^{[n]}(F)/H = Z^n_{eff}(X)/H$.
\end{Cor}
\begin{proof}
This follows immediately from lemma \ref{hilb_chow}.
\end{proof}

It is useful to have a relative version of lemma \ref{hilb_chow} for
flat cycles over a curve:
\begin{Lem} \label{hilb_chow_curve}
Suppose $X/F$ is a projective variety, $C/F$ a curve and $\alpha \in
C_1(X \times C)$ is an effective cycle such that every component of
the support of $\alpha$ is finite and flat over $C$. Then $\alpha =
[Z]$ for some subscheme $Z \subset X \times C$ with $Z \ra C$ flat.
\end{Lem}
Note that by the universal property of the Hilbert scheme, this $Z$
must come from a morphism $C \ra X^{[n]}$ by pulling back the
universal family. 
\begin{proof}
Consider the restriction $\alpha'$ of the cycle $\alpha$ to $X_{F(C)}$
(the generic fiber of the family $X \times C$). We have $\alpha' \in
Z^n_{eff}(X_{F(C)})$ for some $n$, and so we may use lemma
\ref{hilb_chow} to find a subscheme $Z' \in X_{F(C)}$ representing
it. We may interpret $Z'$ as a point in $X_{F(C)}^{[n]}(F(C)) =
X^{[n]}(F(C))$, and therefore obtain a rational morphism $Spec(F(C))
\ra X^{[n]}$. Since $C$ is a curve and $X^{[n]}$ is proper, we may
complete this to a morphism $C \ra X^{[n]}$, and hence obtain a family
$Z \subset X \times C$.

By construction it is clear that $[Z]$ and $\alpha$ both have the same
restriction to $X_{F(C)}$.  We may therefore find an open subset $U
\subset C$ such that the cycles $\alpha$ and $[Z]$ are equal. From the
fundamental sequence
$$C_1(X \times (C \setminus U)) \ra C_1(X \times C) \ra C_1(X \times
U) \ra 0,$$ we see that the difference cycle $[Z] - \alpha$ is supported
entirely on $C_1(X \times (C \setminus U))$. But since the support of
each cycle is flat over $\bP^1$, there cannot be any components
supported in over $\bP^1 \setminus U$, and therefore $[Z] - \alpha =
0$ as claimed.
\end{proof}

If $X/F$ is a projective variety and $L/F$ a finite field extension of
degree $n$. We define a map of sets
\begin{align*}
\cH : X(L) &\ra Z^n_{eff}(X) \\
(\phi : Spec(L) \ra X) &\mapsto \phi_*[Spec(L)]
\end{align*}

\begin{Lem} \label{new_transfer}
Let $X/F$ be a projective variety, $L/F$ a finite field extension, and
suppose we have $x, y \in X^{[n]}(L)$ with $x \sim_H y$. Then $\cH(x)
\sim_H \cH(y)$. In the case $x, y \in X^{(n)}(L)$ are elementarily
linked, so are $\cH(x)$ and $\cH(y)$.
\end{Lem}
\begin{proof}
It suffices to consider the case where $x$ and $y$ are elementarily
H-linked. Therefore, we may reduce either to the case that $x \sim_R
y$ or $[x] = [y]$. If $[x] = [y]$, we may write $[x] = [y] = n[z]$ for
$z \cong Spec(E)$ an irreducible subscheme, and $E \subset L$ a
subfield, $n = [L:E]$. Therefore $x$ and $y$ may only differ by an
element of $Gal(L/F)$ and so $\cH(x) = \cH(y)$.

We may therefore assume that $x$ and $y$ are elementarily
linked. Choose $\phi : \bP^1_E \ra X$ with $\phi(0) = x, \phi(\infty)
= y$, and let $\rho : \bP^1_E \ra \bP^1$ be the natural
covering. Since the cycle $(\phi \times \rho)_*[\bP^1_E] \in C_1(X
\times \bP^1)$ satisfies the conditions of lemma
\ref{hilb_chow_curve}, it follows from lemma \ref{hilb_chow_curve} and
the remark just following it that we can find a morphism $\psi: \bP^1
\ra X^{[n]}$, where $n = [E:F]$ such that if $C \subset \bP^1 \times
X$ is the corresponding family, $[C] = (\phi \times \rho)_*[\bP^1_E]$.

If we denote by $i_p : Spec(F) \ra \bP^1$, $p = 0, \infty$ the
inclusion of points on $\bP^1$, and consider the pullback diagram:
\begin{equation}
\xymatrix{
Spec(E) \ar[r] \ar[d]_x & \bP^1_E \ar[d]^{\phi \times \rho} \\
X \ar[r] \ar[d] & X \times \bP^1 \ar[d] \\
Spec(F) \ar[r]^{i_0} & \bP^1.}
\end{equation}

We have $\cH(x) = x_*(Spec(E)) = x_* i_0^![\bP^1_E]$ which may be
rewritten using \cite{Ful:IT}, theorem 6.2 as $i_0^!(\phi \times
\rho)_*[\bP^1_E] = i_0^![C]$ which by \cite{Ful:IT}, section 10.1 is
the same as $[i_0^{-1}(C)] = [\psi(0)]$, and similarly $\cH(y) =
[\psi(\infty)]$, showing that these points are elementarily H-linked.
\end{proof}

Suppose $X/F$ is a projective variety. Given a zero-dimensional
subscheme $i : z \hra X^{[n]}$, we obtain a family $\cF \subset z \times
X$. We define the cycle $[n]_*(z) \in Z(X)$ by the formula $[n]_*(z) =
{\pi_2}_*[\cF]$.

\begin{Lem} \label{iterated_H}
Let $X/F$ be a projective variety. Then the map $X^{[n][m]}(F) \ra
Z^{nm}_{eff}(X)$ defined by mapping a degree $m$ scheme $z \subset
X^{[n]}$ to $[n]_*[z]$ passes to H-equivalence.
\end{Lem}
\begin{proof}
To show this, it suffices to show that if we have $\phi : \bP^1 \ra
X^{[n][m]}$, $\phi(0) = z$, $\phi(1) = z'$ then $[n]_*[z] \sim_H
[n]_*[z']$. To see this, we will construct a morphism $\psi : \bP^1
\ra X^{[mn]}$ such that $[\psi(0)] = [n]_*[z]$ and $[\psi(\infty)] =
[n]_*[z']$. By the universal property of the Hilbert scheme, this
means that we really need to construct a family $\til{W} \subset X
\times \bP^1$ whose specializations over $0$ and $\infty$ are
$[n]_*[z]$ and $[n]_*[z']$ respectively.

Consider the family corresponding to the map $\phi$. This is a
subscheme $Z \subset X^{[n]} \times \bP^1$ with fibers $z$ and $z'$
over the points $0$ and $\infty$ respectively. Pulling back the
universal family on $X^{[n]}$ via the morphism $Z \ra X^{[n]}$, we
obtain a family $W \hra X \times \bP^1 \times Z$, which is degree $mn$
over $\bP^1$, and such that each component of $W$ is flat over
$\bP^1$. By lemma \ref{hilb_chow_curve}, we may find $\til{W} \subset
X \times \bP^1$ such that $[\til{W}] = \pi_*[W]$, where $\pi : X
\times \bP^1 \times Z \ra X \times \bP^1$ is the projection. It is now
routine to check that the fibers over $0$ and $\infty$ of $\til{W}$
give subschemes whose cycles are equal to $[n]_*[z]$ and $[n]_*[z']$
respectively.
\end{proof}

\begin{Def}
For a scheme $X$ and positive integers $n, m$, we define $X^{(n, m)}$
to be the fiber product:
$$
\xymatrix{
X^{(n, m)} \ar[r] \ar[d]_\pi & S^m \left(S^{n} X\right) \ar[d] \\
X^{(nm)} \ar[r] & S^{nm} X}
$$
\end{Def}

\begin{Lem} \label{connect_lem_new}
Suppose $F$ is prime to $p$ closed, and $X/F$ a quasiprojective
variety. Then the natural morphism $\pi: X^{(n, m)} \ra X^{(nm)}$ is
surjective on $F$-points whenever $n, m$ are powers of $p$.
\end{Lem}
\begin{proof}
Since we may identify $S^m S^{n} X$ with the quotient
$$(X^{nm})/\big((S_{n})^{m} \rtimes S_{m}\big),$$ it follows that the
degree of the map $\pi$ is $\frac{nm!}{(n!)^{m}(m!)}$ which is prime
to $p$ (recall that $v_p(p^r!) = \frac{p^r - 1}{p-1}$ where $v_p$ is
the $p$-adic valuation). Since $\pi$ factors through the \'etale map
$X^{m+n}_\circ \ra X^{(nm)}$ it is also \'etale. In particular, if $x \in
X^{(nm)}(F)$, the fiber $\pi^{-1}(x)$ is \'etale over $Spec(F)$ and
hence the spectrum of a direct sum of separable field extensions
$\oplus L_i$. Since the total degree of this extension is prime to
$p$, there must be at least one of the field extensions $L_i$ whose
degree is not a multiple of $p$. But since $F$ is prime to $p$ closed,
this implies that $L_i = F$, and so the fiber has an $F$-point as
desired.
\end{proof}

\begin{Cor} \label{iterative_cor}
Let $X/F$ be a projective variety.  There is a natural map
$X^{[n][m]}(F)/H \ra X^{[nm]}/H$. In the case $ind(X) = mn$, we have a
natural map $X^{(n)(m)}(F)/H \ra X^{(mn)}/H$. If we also have that $F$
is prime to $p$ closed and $m,n$ are $p$-powers then the map
$X^{(n)(m)}(F)/H \ra X^{(mn)}/H$ is surjective.
\end{Cor}
\begin{proof}
This is immediate from lemmas \ref{index_hilb_open}, \ref{iterated_H},
\ref{connect_lem_new} and corollary \ref{cycle_hilb_H}.
\end{proof} 

\begin{Lem} \label{alg_h_triv}
Let $F$ be prime to $p$ closed, and suppose $X/F$ is a projective
variety. Fix a $p$-power $n$. Suppose $X_L^{(n)}$ is H-trivial for
every finite field extension $L/F$ of $p$-power degree $m$. Then
$X^{(nm)}$ is H-trivial. In particular, we show that if $\alpha \in
X^{(n)}(F)$, then for all $\beta \in X^{(mn)}(F)$, we have $\beta
\sim_H m[\alpha]$.
\end{Lem}
\begin{proof}
By corollary \ref{iterative_cor}, it is sufficient to show that
$X^{(n)(m)}$ is H-trivial. Choose $\alpha \in X^{(n)}(F)$, which is
nonempty by the hypothesis. We will show that given $\beta \in
X^{(n)(m)}(F)$, we can write $[\beta] \sim_H m[\alpha]$. We may write
$\beta = \cH(\til{\beta})$ for some $\til{\beta} \in X^{(n)}(E)$ where
$E/F$ is a degree $m$ \'etale extension.  Choose $\alpha \in
X^{(n)}(F)$, and define $\til{\alpha} \in X^{(n)}(E)$ via composing
$\alpha$ with the structure morphism $Spec(E) \ra Spec(F)$. We then
have $\cH(\til{\alpha}) = n[\alpha]$. Since $X_E^{(n)}$ is H-trivial,
$\til{\alpha} \sim_H \til{\beta}$ and by lemma \ref{new_transfer},
$m[\alpha] \sim_H [\beta]$ as desired.
\end{proof}

\begin{Cor} \label{ind_triv}
Suppose $X/F$ is a projective variety with $F$ is prime to $p$-closed
and such that for every finite field extension $L/F$, $X_L^{(ind
X_L)}$ is H-trivial. Then for every $p$-power $n \geq ind(X)$,
$X^{(n)}$ is H-trivial.
\end{Cor}
\begin{proof}
By lemma \ref{alg_h_triv}, it suffices to show that $X_L^{ind(X)}$ is
H-trivial, where $[L:F] = n/ind(X)$. We prove this by induction on
$ind(X)$. If $ind(X) = 1$, the hypothesis implies that $X_E$ is
R-trivial for every extension $E/F$ and the conclusion follows from
lemma \ref{alg_h_triv} (setting $n = 1$ in the statement of the
lemma).

For the general induction case, we either have $ind(X_L) = ind(X)$ or
$ind(X_L) < ind(X)$. In the first case, the hypothesis immediately
implies $X_L^{ind(X)} = X_L^{ind(X_L)}$ is H-trivial. In the latter
case, we have $m(ind(X_L)) = ind(X)$ for some $p$-power $m$, and by
lemma \ref{alg_h_triv}, to show that $X_L^{ind(X)}$ is H-trivial, it
suffices to show that $X_E^{ind(X_L)}$ is H-trivial for $E/L$ a degree
$m$ extension. Therefore, the result follows from the induction step.
\end{proof}

\begin{Thm} \label{ch_thm}
Suppose $X/F$ is a projective variety with $F$ is prime to $p$ closed,
$p \neq char(F)$ or $F$ perfect and such that $(X_L)^{(ind(X_L))}$
H-trivial for every finite field extension $L/F$. Then $A_0(X) = 0$.
\end{Thm}
\begin{proof}
Let $\alpha \in X^{(i)}$, where $i = ind(X)$. Since by assumption on
the characteristic every prime cycle $\beta$ is represented by a point
in $X^{(n)}(F)$ for some $p$-power $n$, it follows from corollary
\ref{ind_triv} and lemma \ref{alg_h_triv} that $[\beta] \sim_H
\frac{n}{i} [\alpha]$. In particular, $CH_0(X) \cong \bZ$, generated
by $[\alpha]$.
\end{proof}

\begin{Def}
Suppose $f: X \ra Y$ is a morphism of $F$-schemes. We say that $f$ has
R-trivial fibers if for every field extension $L/F$ and every point $y
\in Y(L)$ the fiber $X_y$ is an R-trivial $L$-scheme. Here $X_y$ is
the scheme-theoretic fiber defined as the pullback of $f$ along the
morphism $y : Spec(L) \ra Y$.
\end{Def}

\begin{Def}
We define an equivalence relation on projective varieties which we
call stable R-isomorphism to be the equivalence relation generated by
setting $X$ and $Y$ to be equivalent if there exists $f : X \ra Y$
with R-trivial fibers.
\end{Def}

\begin{Prop} \label{rtriv_proj_fibers}
Suppose $f : X \ra Y$ is an morphism between quasiprojective varieties
with R-trivial fibers. Then $ind(X) = ind(Y)$. If we let $m = ind(X) =
ind(Y)$, then there is an induced set map $X^{(m)}(F) \ra Y^{(m)}(F)$
which is surjective on R-equivalence classes and injective in the case
$X, Y$ are projective. In particular, if $m = 1$, and $X$ and $Y$ are
projective, $f$ is bijective on $R$-equivalence classes.
\end{Prop}

\begin{Cor} \label{stably_rtriv}
If $X$ and $Y$ are stably R-isomorphic then $ind(X) = ind(Y)$ and
there is a bijection $X^{(m)}(F)/R \leftrightarrow Y^{(m)}(F)/R$ where
$m = ind(X) = ind(Y)$.
\end{Cor}

\begin{proof}[proof of proposition \ref{rtriv_proj_fibers}]
It is clear since there is a morphism from $X$ to $Y$ that $ind(Y)$ divides
$ind(X)$. Since the fibers of $f$ are nonempty, we also have that
$ind(X) | ind(Y)$. 

Let $x \in X^{(m)}(F)$. Considering $x \subset X$ as a finite
subscheme, we see as in the proof of lemma \ref{index_hilb_open}, that
$x = Spec(L)$ for some field extension $L/F$ of degree $m$. If we
consider the image $f(x)$, we find similarly that $f(x) \cong x$,
since otherwise the image would have smaller degree, contradicting
$ind(X) = ind(Y)$. Therefore, $f$ induces a map $X^{(m)} \ra
Y^{(m)}$. We note that this may also be seen as the rational map
induced by the morphism $S^m X \ra S^m Y$. Since every $x \in X^{(m)}$
is of the form $Spec(L)$ for $L/F$ a degree $m$ field extension, we
have a commutative diagram such that the vertical arrows are
surjective:
\begin{equation*}
\xymatrix{ \underset{[L:F]}{\coprod} X(L) \ar[d]_{\cH}
\ar[rr]^{\coprod f_L} & & \underset{[L:F]}{\coprod} Y(L) \ar[d]^{\cH}
\\ X^{(m)}(F) \ar[rr] & & Y^{(m)}(F) }
\end{equation*}
It is clear by tracing the diagram that the map on the bottom must be
surjective, and it is also clear that it must preserve R-equivalence
classes. 

We need only show therefore that the map is injective on R-equivalence
classes in the case $X$ and $Y$ are projective. Without loss of
generality, we may assume that we have $x, x' \in X^{(m)}(F)$ such
that $y = f(x)$ and $y' = f(x')$ are elementarily linked. Choose
$\bP^1 \ra Y^{(m)}$ linking $y$ and $y'$. In the case that $y = y'
\cong Spec(L)$, we have that $x, x'$ may be lifted to elements of
$X(L)$ which both lie in the same fiber over a point in $Y(L)$. Since
by hypothesis, the fibers of $f$ are R-trivial, we therefore have $x
\sim_R x'$ due to the fact that $\cH$ preserves R-equivalence in this
case.

We are therefore done if we may show that there is a morphism $\bP^1
\ra X^{[m]}$ connecting some point in the fiber over $y$ with a point
over the fiber of $y'$. We begin by choosing a morphism $\phi: \bP^1
\ra Y^{[m]}$ connecting $y$ and $y'$, and consider the pullback of the
universal family. This gives a curve $C \subset Y \times \bP^1$ such
that the projection $C \ra \bP^1$ is degree $m$ and such that the
fibers over $0$ and $\infty$ are equal to $y$ and $y'$
respectively. If we consider the projection morphism $C \ra Y$
restricted to the generic point $Spec(F(C))$, we obtain a point in
$Y(F(C))$. Since $f$ has R-trivial fibers, the fiber over this point
in $X$ is nonempty and hence there is a morphism $Spec(F(C)) \ra X$
such that its composition with $f$ gives the original map $Spec(F(C))
\ra Y$. Let $\til{C} \ra C$ be the normalization of $C$. Since $X$ is
projective, we get a morphism $\til{C} \ra X$ such that the digram
\begin{equation*}
\xymatrix{
\til{C} \ar[r] \ar[dr] & X \times \bP^1 \ar[r] \ar[d] & Y \times \bP^1
\ar[dl] \\
 & \bP^1}
\end{equation*}
is commutative. In particular, if we let $D \subset X \times \bP^1$ be
the image of $\til{C}$, then $D$ is birational to $C$ and by the
universal property of $X{[m]}$, defines a morphism $\bP^1 \ra
X{[m]}$. If we let $U \subset \bP^1$ be the open set on which $C \ra
\bP^1$ is \'etale, one may check that we have a commutative diagram

\begin{equation*}
\xymatrix{
U \ar[r]^{\phi|_U} \ar[rd] & X^{(m)} \ar@{.>}[d] \\
& Y^{(m)}}
\end{equation*}

Which shows that we may connect points in the fiber over $y$ and $y'$ as desired.
\end{proof}

\begin{Lem} \label{unirational_fibers}
Suppose $f : X \ra \bP^1$ is a dominant morphism of quasiprojective
varieties such that the fibers are unirational of constant positive
dimension. Then we may find $x, y \in X(F)$ such that $f(x) = 0, f(y)
= \infty$ and $x$ and $y$ are elementarily linked.
\end{Lem}
\begin{proof}
Since the generic fiber is unirational, we find a rational map $g :
\bP^N \times \bP^1 \ra X$ commuting with $f$. By the hypothesis, we
may assume that $g$ is defined on an open subset of codimension at
least $2$, and consequently, at some point in each fiber over
$\bP^1$. Choose points $x', y' \in \bP^N(F)$ such that $g$ is defined
at $(x', 0)$ and $(y', \infty)$. Choosing a linear map $\phi : \bP^1
\ra \bP^N$ with $\phi(0) = x'$ and $\phi(\infty) = y'$, we find $g
\circ \phi$ elementarily links the points $x = g(x', 0)$ and $y =
g(y', \infty)$ as desired.
\end{proof}

\begin{Cor} \label{unirat_fibers}
Suppose $f : X \ra Y$ is a morphism with R-trivial fibers which are
unirational of constant positive dimension. Then $f$ is bijective on
R-equivalence classes.
\end{Cor}
\begin{proof}
Since $f$ is clearly surjective on R-equivalence classes, we need only
to show that it is injective. It suffices to consider the case that
$x, x' \in X$ with $f(x)$ and $f(x')$ elementarily linked. But
considering a path $\phi : \bP^1 \ra Y$ linking the two points, we may
pullback the family $f: X \ra Y$ over $\phi$ to obtain a family over
$\bP^1$. Hence by lemma \ref{unirational_fibers}, we can reduce to the
case that $f(x) = f(x')$. But in this case we are done since the
fibers of $f$ are R-trivial by hypothesis. 
\end{proof}

\section{Preliminaries on Severi-Brauer flag varieties} \label{sbv_flag}

\begin{Def}
Let $A$ be an central simple algebra of degree $n$. Choose positive
integers $n_1 < \ldots < n_k < n$. The Severi-Brauer flag variety of
type $(n_1, \ldots, n_k)$, denoted $V_{n_1, \ldots, n_k}(A)$, is the
variety whose points correspond to flags of ideals $I_{n_1} \subset
I_{n_2} \subset \cdots \subset I_{n_k}$, where $I_{n_i}$ has reduced
dimension $n_k$. More precisely $V_{n_1, \ldots, n_k}(A)$ represents
the following functor:
$$V_{n_1, \ldots, n_k}(A)(R) = \left\{(I_1, \ldots, I_k) 
\left|
\begin{matrix}
I_i \in Gr(n_i n, A)(R) \\
\text{ is a right ideal of $A_R$ and $I_i \subset I_{i+1}$}
\end{matrix}
\right. \right\}
$$
\end{Def}

In particular, in the case $k = 1$, the variety $V_i(A)$ is the $i$'th
generalized Severi-Brauer variety of $A$ (\cite{Blanchet}), which
parametrizes right ideals of $A$ which are locally direct summands of
reduced rank $i$. The same definition generalizes easily to sheaves of
Azumaya algebras of constant degree over a base scheme $S$.

\begin{Thm} \label{flag_reduction}
Suppose $A$ is a central simple $F$-algebra. Then the Severi-Brauer
flag variety $V_{n_1, \ldots, n_k}(A)$ is stably R-isomorphic to
$V_d(D)$, where $D$ is any central simple algebra Brauer equivalent to
$A$ and
$$d = GCD\{n_1, \ldots, n_k, ind(A)\}.$$
\end{Thm}

This result relies on a number of intermediate results:

\begin{Prop} \label{flag_rtriv}
Suppose $A$ is a central simple $F$-algebra, and we have positive
integers $n_1 < \cdots < n_k$ such that $ind(A) | n_i$, each $i$. Then
any two points on the Severi-Brauer flag variety $V_{n_1, \ldots,
n_k}(A)$ are elementarily linked. In particular, $V_{n_1, \ldots,
n_k}(A)$ is R-trivial.
\end{Prop}

\begin{Cor} \label{gsbv_rtriv}
Suppose $A$ is a central simple $F$-algebra, and we have positive
integers $n$ such that $ind(A) | n$. Then any two $F$-points in the
generalized Severi-Brauer variety $V_n(A)$ are elementarily linked. In
particular, $V_n(A)$ is R-trivial.
\end{Cor}

\begin{Lem} \label{ideal_lem}
Suppose $A = End_{r.D}(V)$ for some $F$-central division algebra $D$,
where $V$ is a right $D$-space. Let $i = deg(D)$, and let $I \subset
A$ be a right ideal of reduced dimension $ri$ (note that every ideal
has reduced dimension a multiple of $i$). Then there exists a
$D$-subspace $W \subset V$ of dimension $r$ such that $I =
Hom_{r. D}(V,W) \subset End_{r.D}(V)$.

Equivalently, writing $A = M_m(D)$, we may consider $I$ to be the set
of matricies such that each column is a vector in $W$.
\end{Lem}
\begin{proof}
Choose a right ideal $I \subset A$, and let $W = im(I)$. It is enough
to show that $I = Hom_{r. D}(V, W)$. The claim concerning reduced
dimension will follow immediately from a dimension count. Since it is
clear by definition that $I \subset Hom_{r. D}(V,W)$, it remains to
show that the reverse inclusion holds. We do this by showing that $I$
contains a basis for $Hom_{r. D}(V, W)$. Let $e_1, \ldots, e_m$ be a
basis for $V$, and $f_1, \ldots, f_r$ be a basis for $W$. We must show
that the transformation $T_{i,j} \in I$ where $T_{i,j}(e_k) = f_j
\delta_{i,k}$. Since $f_i \in im(i)$ for some $i \in I$, we know that
there exists $v \in V$ such that $i(v) = f_j$. Now define $a \in
End_{r.D}(V)$ to be given by $a(e_k) = v \delta_{i,k}$. Now, $ia(e_k)
= f_j \delta_{i,k}$ and $ia \in I$ as desired.
\end{proof}

\begin{proof} [proof of \ref{flag_rtriv}]
Let $A = M_m(D)$ for some division algebra $D$ with $deg(D) = i =
ind(A)$, $mi = n = deg(A)$. We may therefore write $A = End_{r. D}(V)$
for some right $D$-space $V$ of dimension $m$. Choose flags of ideals
$$(I_1, \ldots, I_k), (I_1', \ldots, I_k') \in V_{n_1, \ldots,
n_k}(A)(F).$$ We will show that there is a rational morphism $f: \bA^1
\dra V_{n_1, \ldots, n_k}(A)$, such that $f(0) = (I_1, \ldots, I_k)$
and $f(1) = (I_1', \ldots, I_k')$.

By lemma \ref{ideal_lem}, we may write $I_j = Hom(V, W_j)$, $I_j' =
Hom(V, W_j')$. Choose bases $w_{j,1}, \ldots, w_{j, l_j}$ for $W_j$
and $w_{j, 1}', \ldots, w_{j, l_j}'$ for $W_j'$, where $l_j = n_j /
i$. Define morphisms $f_{j,l} : \bA^1 \ra V$ by $f_{j,l}(t) = w_{j,l} t
+ w_{j,l}' (1-t)$. We may combine these to get rational morphisms $\bA^1
\dra Gr(n_j n, A)$ by taking $t$ to the $n_j n$-dimensional space of
matricies in $M_m(D)$ whose columns are right $D$-linear combinations
of the vectors $w_{j, 1} t + w_{j, 1}' (1 - t), \ldots, w_{j, l_j} t +
w_{j, l_j}' (1 - t)$. By \ref{ideal_lem}, this corresponds to a
rational morphism $f_j : \bA^1 \ra V_{n_j}(A)$. One may check that
$f_j(0) = I_j$ and $f_j(1) = I_j'$. Further, for any $t$, $f_j(t)
\subset f_{j+1}(t)$. Therefore, we may put these together to yield a
rational morphism $f : \bA^1 \dra V_{n_1, \ldots, n_k}(A)$ with $f(0) =
(I_1, \ldots, I_k)$ and $f(1) = (I_1, \ldots, I_k)$.
\end{proof}

\begin{Rem} \label{rtriv_flag_fibers}
In fact the proof above shows that if we are given $n < m$ with
$ind(A) | n, m$, and we fix $I \in V_n(A)$ $\pi^{-1}(I)$ is R-trivial
where $\pi : V_{n,m}(A) \ra V_n(A)$ is the natural projection. In
fact, given $\alpha, \beta \in \pi^{-1}(I)$, the path constructed in
the proof above to connect $\alpha$ and $\beta$ as points in
$V_{n,m}(A)$ lies entirely in the fiber $\pi^{-1}(I)$ showing they are
R-equivalent there as well.
\end{Rem}

\begin{proof}[proof of theorem \ref{flag_reduction}]
Let $X = V_{n_1, \ldots, n_k}(A)$ and $Y = V_d(D)$. Consider the
product variety $X \times Y$ together with its natural projections
$\pi_1, \pi_2$ onto $X$ and $Y$ respectively. I claim that both
projections have R-trivial fibers, which would prove the theorem.

Suppose we have $x : Spec(L) \ra X$ or $x : Spec(L) \ra Y$. This would
imply that $X(L) \neq \emptyset$ or $Y(L) \neq \emptyset$, and in
either case this in turn says that $ind(A) | d$. Since the scheme
theoretic fiber over $x$ is isomorphic to either $X_L$ or $Y_L$
respectively, we know that since $ind(A_L) = ind(D_L) | d$ that the
fibers are R-trivial by proposition \ref{flag_rtriv}.
\end{proof}

\begin{Def}
Suppose $A$ is a central simple algebra and $I \subset A$ is a right
ideal of reduced dimension $l$. Given integers $n_1, \ldots, n_k < l$,
we define the variety $V_{n_1, \ldots, n_k}(I)$ to be the subvariety
of $V_{n_1, \ldots, n_k}(A)$ consisting of flags of ideals all of
which are contained within $I$.
\end{Def}

For these varieties, we have a theorem which generalizes a result from
\cite{Ar:BS} on Severi-Brauer varieties:
\begin{Thm} \label{sub_flag}
Suppose $A$ is a central simple algebra. Let $I \subset A$ be a right
ideal of reduced dimension $l$. Then there exists a degree $l$ algebra
$D$ which is Brauer equivalent to $A$ such that for any $n_1, \ldots,
n_k < k$,
\begin{equation*}
V_{n_1, \ldots, n_k}(I) = V_{n_1, \ldots, n_k}(D)
\end{equation*}
\end{Thm}

In order to prove this theorem, we will use the following lemma:

\begin{Lem} \label{idempotent}
Let $A$ be an Azumaya algebra with center $R$, a Noetherian
commutative ring, and suppose that $I$ is a right ideal of $A$ such
that $A/I$ is a projective $R$-module. Then there exists an idempotent
element $e \in I$ such that $I = eA$.
\end{Lem}
\begin{proof}
Since $A/I$ is projective as an $R$-module, by \cite{DeIn} it is also
a projective (right) $A$-module. This implies that the short exact
sequence
$$0 \ra I \ra A \ra A/I \ra 0$$ splits as a sequence of right
$A$-modules, and therefore, there exists a right ideal $J \subset A$
such that $A = I \oplus J$. We may therefore uniquely write $1 = e +
f$, with $e \in I$ and $f \in J$. Now,
\begin{align*}
e = (e + f)e = e^2 + fe.
\end{align*}
Since $e^2 \in I, fe \in J$ this gives $fe \in I \cap J = 0$ and so
$e^2 = e$. Finally, $I = (e + f)I = eI + fI$, and this gives
$fI \in J \cap I = 0$. Consequently, we have $eA \subset I = eI
\subset eA$ so $I = eA$ as desired.
\end{proof}

\begin{proof}[proof of theorem \ref{sub_flag}]
By \ref{idempotent}, we know that $I = eA$ for some idempotent $e \in
A$. Set $D = eAe$. 

Let $X_I = V_{n_1, \ldots, n_k}(I)$, and $X_D = V_{n_1, \ldots,
n_k}(D)$. To prove the theorem, we will construct mutually inverse
maps (natural transformations of functors) $\phi : X_I \ra X_D$ and
$\psi : X_D \ra X_I$. For a commutative Noetherian $F$-algebra $R$,
and for $\cJ = (J_1, \ldots, J_k) \in X_I(R)$, we define $\phi(\cJ) =
(J_1e, \ldots, J_ke) = (eJ_1e, \ldots, eJ_ke)$. For $\cK = (K_1,
\ldots, K_k) \in X_D(R)$, define $\phi(\cK) = (K_1A_R, \ldots, K_k
A_R)$. To see that these are mutually inverse, we need to show that
for each $i = 1, \ldots, k$, we have $J_ieA = J_i$ and that $K_iA_Re =
K_i$. For the second we have
$$K_iA_Re = K_ie = K_i$$ since $K_i\subset eA_Re$. For the first, we
note that by the lemma, we have $J_i = hA_i$ for some idempotent
$h$. But then
$$J_i \supset J_ieA_R = J_iI \supset J_i^2 = hA_RhA_R = hA_R = J_i$$
and so $J_i = J_ieA_R$ and we are done.
\end{proof}

\section{Moduli spaces of \'etale subalgebras} \label{moduli}

Let $S$ be a Noetherian scheme, and let $A$ be a sheaf of Azumaya
algebras over $S$. Our goal in this section is to study the functor
$\et(A)$, which associates to every $S$-scheme $X$, the set of sheaves
of commutative \'etale subalgebras of $A_X$. We will show that this
functor is representable by a scheme which may be described in terms
of the generalized Severi-Brauer variety of $A$.

Unless said otherwise, all products are fiber products over $S$. If
$X$ is an $S$-scheme with structure morphism $f : X \ra S$, then we
write $A_X$ for the sheaf of $\mathcal{O}_X$-algebras $f^*(A)$. For a
$S$-scheme $Y$, we occasionally write $Y_X$ for $Y \times X$, thought
of as an $X$-scheme.

Every sheaf of \'etale subalgebras may be assigned a discrete
invariant, which we call its type, and therefore our moduli scheme is
actually a disjoint union of other moduli spaces.

To begin, let us define the notion of type. 

\begin{Def}
Let $R$ be a local ring, and $B/R$ an Azumaya algebra. If $e \in B$ is
an idempotent, we define the rank of $e$, denoted $r(e)$ to be the
reduced rank of the right ideal $eB$.
\end{Def}

Let $E$ be a sheaf of \'etale subalgebras of $A/S$, and let $p \in
S$. Let $R$ be the local ring of $p$ in the \'etale topology (so that
$R$ is a strictly Henselian local ring). Then taking \'etale stalks,
we see that $E_p$ is an \'etale subalgebra of $A_p/R$, and it follows
that
$$E_p = \oplus_{i=1}^k R e_i,$$ for a uniquely defined collection of
idempotents $e_i$, which are each minimal idempotents in $S_p$.

\begin{Def}
The type of $E$ at the point $p$ is the unordered collection of
positive integers $[r(e_1), \ldots, r(e_m)]$.
\end{Def}

\begin{Def}
We say that $E$ has type $[n_1, \ldots, n_m]$ if if has this type for
each point $p \in S$.
\end{Def}

\begin{Rem} \label{idempotents}
Since $1 = \sum e_i$, the ideals $I_i = e_i A$ span $A$. Further it is
easy to see that the ideals $I_i$ are linearly independent since $e_ia
= e_jb$ implies $e_i a = e_i e_i a = e_i e_j b = 0$. We therefore know
that the numbers making up the type of $E$ give a partition of
$deg(A_p)$.
\end{Rem}

Some additional notation for partitions will be useful. Let $\rho =
[n_1, \ldots, n_m]$. For a positive integer $i$, let $\rho(i)$ be the
number of occurrences of $i$ in $\rho$. Let $S(\rho)$ be the set of
distinct integers $n_i$ occurring in $\rho$, and let $N(\rho) =
|S(\rho)|$. Let
$$\ell(\rho) = \underset{i \in S(\rho)}{\sum} \rho(i) = m$$
be the length of the partition.

Suppose $A/S$ is an sheaf of Azumaya algebras, and suppose $S$ is a
connected, Noetherian scheme.  Let $\rho = [n_1, \ldots, n_m]$ be a
partition of $n = deg(A)$. Let $\et_{\rho}(A)$ be the functor which
associates to every $S$ scheme $X$ the set of \'etale subalgebras of
$A_X$ of type $\rho$. That is, if $X$ has structure map $f: X \ra S$, 
$$\et_{\rho}(A)(X) = 
\left\{
\begin{matrix}
\text{sub-$\mathcal{O}_X$-modules } \\
E \subset f^*A 
\end{matrix}
\left|
\begin{matrix}
\text{$E$ is a sheaf of commutative \'etale} \\
\text{subalgebras of $f^*A$ of type $\rho$}
\end{matrix}
\right. \right\}$$

Our first goal will be to describe the scheme which represents this
functor. We use the following notation: 
$$V(A)^{\rho} = \underset{i \in S(\rho)}{\prod} V_i(A)^{\rho(i)}.$$

We define $V(A)^{\rho}_*$ to be the open subscheme parametrizing
ideals which are linearly independent. That is to say, for a
$S$-scheme $X$, if $I_1, \ldots, I_{\ell{\rho}}$ is a collection of
sheaves of ideals in $A_X$, representing a point in $V(A)^{\rho}(X)$,
then by definition, this point lies in $V(A)^{\rho}_*$ if and only if
$\oplus I_i = A$.

Let $S_\rho$ be the subgroup $\prod_{i \in S(\rho)} S_{\rho(i)}$ of
the symmetric group $S_n$. For each $i$, we have an action of
$S_{\rho(i)}$ on $V_i(A)^{\rho(i)}$ by permuting the factors. This
induces an action of $S_{\rho}$ on $V(A)^{\rho}$, and on
$V(A)^{\rho}_*$. Denote the quotients of these actions by
$S^{\rho}V(A)$ and $V(A)^{(\rho)}_*$ respectively. We note that since
the action on $V(A)^{\rho}_*$ is free, the quotient morphism
$V(A)^{\rho}_* \ra $ is a Galois covering with group
$S_\rho$.

\begin{Thm}
Let $\rho = [n_1, \ldots, n_m]$ be a partition of $n$. Then the
functor $\et(A)_\rho$ is represented by the scheme $V(A)^{(\rho)}_*$.
\end{Thm}
\begin{proof}
To begin, we first note that both $\et(A)_\rho$ and the functor
represented by $V(A)^{(\rho)}_*$ are sheaves in the \'etale
topology. Therefore, to show that these functors are naturally
isomorphic, it suffices to construct a natural transformation $\psi :
V(A)^{(\rho)}_* \ra \et(A)_\rho$, and then show that this morphism
induces isomorphisms on the level of stalks.

Let $X$ be an $S$-scheme, and let $p : X \ra V(A)^{(\rho)}_*$. To
define $\psi(X)(p)$, since both functors are \'etale sheaves, it
suffices to define it on an \'etale cover of $X$. Let $\til{X}$ be the
pullback in the diagram
\begin{equation} \label{pullback}
\xymatrix{ \til{X} \ar[r] \ar[d] & V(A)^{\rho}_* \ar[d]^{\pi} \\ X
\ar[r]^p & V(A)^{(\rho)}_*}
\end{equation}

Since the quotient morphism $\pi$ is \'etale, so is the morphism
$\til{X} \ra X$. Therefore we see that after passing to an \'etale
cover, and replacing $X$ by $\til{X}$, we may assume that $p = \pi(q)$
for some $q \in V(A)^{\rho}_*(X)$. Passing to another cover, we may
also assume that $X = Spec(R)$.

Since $p = \pi(q)$, we may find right ideals $I_1, \ldots,
I_{\ell(\rho)}$ of $A_R$ such that $\oplus I_i = A_R$, which represent
$q$. Writing
$$1 = \sum e_i, \ \ e_i \in I_i,$$ we define $E_p = \oplus e_i
R$. This is a split \'etale extension of $R$, which is a subalgebra of
$A$, and we set $\psi(p) = E_p$. One may check that this defines a
morphism of sheaves. Note that this definition with respect to an
\'etale cover gives a general definition since the association $(I_1,
\ldots, I_{\ell}) \mapsto E_p$ is $S_\rho$ invariant.

To see that $\psi$ is an isomorphism, it suffices to check that it is
an isomorphism on \'etale stalks. In other words, we may restrict to
the case that $X = Spec(R)$, where $R$ is a strictly Henselian local
ring. 

We first show that $\psi$ is injective. Suppose $E$ is an \'etale
subalgebra of $A_R$ of type $\rho$. Since $R$ is strictly Henselian,
we have 
$$E = \underset{i \in S(\rho)}{\oplus} \underset{j =
1}{\overset{\rho(i)}{\oplus}} e_{i,j} R.$$ By definition, since the
type of $E$ is $\rho$, if we let $I_{i,j} = e_{i,j}A$, then we the
tuple of ideals $(I_{i,j})$ defines a point $q \in
V(A)^{\rho}_*(R)$. Further, since $\sum e_{i,j} = 1$, we actually have
$q \in V(A)^{\rho}_*(R)$. If we let $p = \pi(q)$, then tracing through
the above map yields $\psi(R)(p) = E$. Therefore $\psi$ is surjective.

To see that it is injective, we suppose that we have a pair of points
$p, p' \in V(A)^{(\rho)}_*(R)$. By forming the pullbacks as in
equation \ref{pullback}, since $R$ is strictly Henselian, we
immediately find that in each case, because $\til{X}$ is an \'etale
cover of $X$, it is a split \'etale extension, and hence we have
sections. This means we may write
$$p = \pi(I_1, \ldots, I_{\ell(\rho)}), \ \ p' = \pi(I_1', \ldots,
I_{\ell(\rho)}').$$ Note that in order to show that $p = p'$ is
suffices to prove that the ideals are equal after reordering. Now, if
$E_p = E_{p'}$, then both rings have the same minimal
idempotents. However, by remark \ref{idempotents}, the ideals are
generated by these idempotents. Therefore, the ideals coincide after
reordering, and we are done.
\end{proof}

Since we now know that the functor $\et_\rho(A)$ is representable, we
will abuse notation slightly and refer to it and the representing
variety by the same name.

\begin{Def}
$\et(A)$ is the disjoint union of the schemes $\et_{\rho}(A)$ as
$\rho$ ranges over all the partitions of $n = deg(A)$.
\end{Def}

\begin{Cor}
The functor which associates to any $S$-scheme $X$ the set of \'etale
subalgebras of $A_X$ is representable by $\et(A)$.
\end{Cor}

\begin{Rem} \label{grass}
By associating to an \'etale subalgebra $E \subset A_X$ its underlying
module, we obtain a natural transformation to the Grassmannian
functor, $\et_\rho(A) \ra Gr(\ell(\rho), A)$. 
\end{Rem}

\section{Subfields of central simple algebras} \label{subfields}

In this section and for the remainder of the paper, we specialize back
to the case where $S = Spec(F)$, and $A$ is a central simple
$F$-algebra. If $E$ is an \'etale subalgebra of $A$, then taking the
\'etale stalk at $Spec(F)$ amounts to extending scalars to the
separable closure $F^{sep}$ of $F$. Let $G$ be the absolute Galois
group of $F^{sep}$ over $F$. Writing
$$E \otimes F^{sep} \cong \underset{i \in S(\rho)}\oplus \underset{j =
1}{\overset{\rho(i)}{\oplus}} e_{i,j} F^{sep},$$ we have an action of
$G$ on the idempotents $e_{i,j}$. One may check that the idempotents
$e_{i,j}$ are permuted by $G$, and there is a correspondence between
the orbits of this action and the idempotents of $E$. In particular we have
\begin{Lem}
In the notation above, if $E$ is a subfield of $A$, then $|S(\rho)| =
1$.
\end{Lem}
\begin{proof}
$E$ is a field if and only if $G$ acts transitively on the set of
idempotents. On the other hand, this action must also preserve the
rank of an idempotent, which implies that all the idempotents have the
same rank.
\end{proof}

Therefore, if we are interested in studying the subfields of a central
simple algebra, we may restrict attention to partitions of the above
type.  If $m | n = deg(A)$, we write
$$\et_m(A) = \et_{[\frac{n}{m}, \frac{n}{m}, \ldots, \frac{n}{m}]}.$$

Note that every separable subfield of dimension $m$ is represented by
a $F$-point of $\et_m(A)$, and in the case that $A$ is a division
algebra, this gives a 1-1 correspondence. In particular, elements of
$\et_n(A)(F)$ are in natural bijection with the maximal separable
subfields of $A$.

\begin{Prop} \label{subfields_unirat}
Suppose $A$ is a central simple $F$ algebra of degree $md = n$, and
suppose $\et_m(A)(F) \neq \emptyset$. Then $\et_m(A)$ is unirational.
\end{Prop}
\begin{proof}
Note that if $F = \ov{F}$, any two \'etale subalgebras of type $[d,
\ldots, d]$ are conjugate under the action of $GL_1(A) =
A^*$. Therefore, if $L \subset A$ is an \'etale subalgebra of the
appropriate type, the morphism $GL_1(A) \ra \et_m(A)$ defined by $g
\ra [gLg^{-1}]$ is dominant. Since $GL_1(A)$ is rational, $\et_m(A)$
is unirational.
\end{proof}

\subsection{Maximal subfields} \label{maxl}

Note that if $a \in A$ is an element whose characteristic polynomial has
distinct roots, then the field $F(a)$ is a maximal \'etale subalgebra
of $A$.

\begin{Thm} \label{max_surject}
Let $U \subset A$ be the Zariski open subset of elements of $A$ whose
characteristic polynomials have distinct roots. Then there is a
dominant rational map $U \ra \et_n(A)$ which is surjective on
$F$-points.
\end{Thm}
\begin{proof}
This argument is a geometric analog of one in \cite{KraSa}.  Consider
the morphism $U \ra Gr(n, A)$ defined by taking an element $a$ to the
$n$-plane spanned by the elements $1, a, a^2, \ldots, a^{n-1}$. Since
the characteristic polynomial of $a$ has distinct roots, this
$n$-plane is a maximal \'etale subalgebra. Therefore by the remark at
the end of section \ref{moduli}, we obtain a morphism $U \ra
\et_n(A)$. This morphism can be described as that which takes an
element of $A$ to the \'etale subalgebra which it generates. Since
every \'etale subalgebra of $A$ can be generated by a single element,
this morphism is surjective on $F$-points. Since this also holds after
fibering with the algebraic closure, it follows also that this
morphism is surjective at the algebraic closure and hence dominant.
\end{proof}

\begin{Thm} \label{maximal}
Suppose $A$ has degree $n$. Then $\et_n(A)$ is R-trivial.
\end{Thm}
\begin{proof}
Since any two points on $A$ as an affine space are elementarily
linked, the open subscheme $U \subset A$ from the previous theorem is
R-trivial. Therefore since $U$ is R-trivial and there is a map $U \ra
\et_n(A)$ which is surjective on $F$-points, it follows that
$\et_n(A)$ is R-trivial as well.
\end{proof}

\subsection{Degree $4$ algebras} \label{deg4}

In this section we assume that $char(F) \neq 2$.

\begin{Lem} \label{degree4_pluker}
Let $A$ be a degree $4$ central simple $F$-algebra. Then $V_2(A)$ is
isomorphic to an involution variety $V(B, \sigma)$ of a degree $6$
algebra with orthogonal involution $\sigma$ (see section \ref{inv} for
the definition of involution varieties).
\end{Lem}
\begin{proof}
Consider the map 
$$Gr(2,4) \ra \mathbb{P}^5,$$ given by the Pl\"uker embedding. Fixing $V$ a
$4$ dimensional vector space, we may consider this as the map which
takes a $2$ dimensional subspace $W \subset V$ to the $1$ dimensional
subspace $\wedge^2W \subset \wedge^2V$. This morphism gives an
isomorphism of $Gr(2,4)$ with a quadric hypersurface. This quadric
hypersurface may be thought of as the quadric associated to the
bilinear form on $\wedge^2V$ defined by $<\omega_1, \omega_2> =
\omega_1 \wedge \omega_2 \in \wedge^4 V \cong F$. Note that one must
choose an isomorphism $\wedge^4 V \cong F$ to obtain a bilinear form,
and so it is only defined up to similarity. Nevertheless, the quadric
hypersurface and associated adjoint (orthogonal) involution depend
only on the similarity class and are hence canonically defined.

Since the Pl\"uker embedding defined above is clearly $PGL(V)$
invariant, using \cite{Ar:BS}, for any degree $4$ algebra $A$ given by a
cocycle $\alpha \in H^1(F, PGL_4)$, we obtain a morphism:
$$V_2(A) \ra V(B),$$ where $B$ is given by composition of $\alpha$
with the standard representation $PGL(V) \ra PGL(V \wedge V)$. By
\cite{Ar:BS} this implies that $B$ is similar to $A^{\otimes2}$ in
$Br(F)$. Also, it is easy to see that the quadric hypersurface and
hence the involution is $PGL_4$ invariant, and hence descends to an
involution $\sigma$ on $B$. We therefore obtain an isomorphism
$V_2(A) \cong V(B, \sigma)$ as claimed.
\end{proof}

\begin{Thm} \label{degree4}
Suppose $A$ is a degree $4$ central simple $F$ algebra. Then
$\et_2(A)$ is R-trivial.
\end{Thm}
\begin{proof}
By lemma \ref{inv_rtriv} (although located near the end of this paper,
it does not require the present result), it follows that any two
points in $V(B, \sigma)^{(2)}(F)$ are elementarily linked. In
particular, since $\et_2(A)(F) = V_2(A)^{(2)}_*(F) \subset
V_2(A)^{(2)}(F) = V(B, \sigma)^{(2)}(F)$ by lemma
\ref{degree4_pluker}, we may conclude that any two points on
$\et_2(A)(F)$ are also elementarily linked.
\end{proof}

\subsection{Exponent $2$ algebras} \label{exp2_subsec}

We assume again in this section that $char(F) \neq 2$. We will show in
this section that for an algebra of exponent $2$, the variety
$\et_{deg(A)/2}(A)$ is R-trivial. It will be useful, however, to prove
the slightly more general fact below:

\begin{Thm} \label{exp2}
Suppose $A$ is an algebra of exponent $2$ and degree $n = 2m$. Then
every nonempty open subvariety $U \subset \et_m(A)$ is R-trivial.
\end{Thm}

The idea of the argument presented in the proof of this theorem is due
to D. Saltman.

For the remainder of the section, fix $A$ as in the hypotheses of the
theorem above. For an involution $\tau$ (symplectic or orthogonal),
let $Sym(A, \tau)$ denote the subspace of elements of $A$ fixed by
$\tau$.

\begin{Lem}
Suppose $a \in A$ generates an \'etale subalgebra $F(a) \in
\et_m(A)$. Then there is an symplectic involution $\sigma$ on $A$ such
that $F(a) \subset A^\sigma$.
\end{Lem}
\begin{proof}
Since $A$ has exponent $2$, we may choose $\tau$ to be an arbitrary
symplectic involution on $A$. Set $L = F(a)$, and let $\phi = \tau |_L
: L \ra A$. By the Noether Skolem theorem (\cite{DeIn}, Cor 6.3),
there is an element $u \in A$ such that conjugation by $u$ restricted
to $F(a)$ gives $\phi$.

I claim we may choose $u$ so that $\tau(u) + u$ is invertible.  If
this is the case, set $v = \tau(u) + u$. Since $u l = \tau(l)u$ for
all $l \in L$, we may take $\tau$ of both sides to obtain $\tau(l)
\tau(u) = \tau(u) l$. Adding opposite sides of these two equations
yields $vl = \tau(l) v$, or in other words $inn_{v^{-1}} \circ \tau
|_L = id_L$. But since $v$ is $\tau$-symmetric, $\sigma = inn_{v^{-1}}
\circ \tau$ is a symplectic involution, and by construction
$\sigma(l) = l$ for $l \in L$, proving the lemma. Hence we need only
prove the claim. This is done as follows:

Suppose $w$ is any element of $A^*$ such that $inn_w |_L = \tau
|_L$. Let $Q = C_A(L)$ be the centralizer of $L$ in $A$. For any $q
\in Q^*$, it is easy to check that $inn_{uq}|_L = inn_u |_L$. Define a
linear map $f : Q \ra A$ by $f(q) = uq + \tau(uq)$. The condition that
$f(q) \in A^*$ is an open condition, defining an open subvariety $U
\in Q$. I claim that $U$ is not the empty subvariety. Note that since
$Q$ is an affine space and $F$ is infinite, then $F$-points are dense
on $Q$ and this would imply that $U$ contains an $F$-point. 

To check that $U$ is not the empty subvariety, it suffices to check
that $U(\ov(F)) \neq \emptyset$. In other words, we must exhibit an
element $q$ in $Q_{\ov{F}} = Q \otimes_F \ov{F}$ such that $uq +
\tau(uq)$ is invertible in $A_{\ov{F}} = A \otimes_F \ov{F}$ (we have
abused notation here by writing $u$ in place of $u \otimes 1$). To do
this, first choose a symplectic involution $\gamma$ on $A_{\ov{F}}$
such that $\gamma|_{\tau(L_{\ov{F}})} = id_{\tau(L_{\ov{F}})}$. By
\cite{BofInv} we may find an element $r \in Symm(A_{\ov{F}}, \tau)$
such that $\gamma \circ \tau = inn_r$. We therefore have
$$inn_r|_{L_{\ov{F}}} = \gamma \circ \tau |_{L_{\ov{F}}} = \tau
|_{L_{\ov{F}}} = inn_w |_{L_{\ov{F}}}.$$ This in turn implies that
$inn_{w^{-1}r}|_{L_{\ov{F}}} = id_{L_{\ov{F}}}$, or in other words
$w^{-1}r \in C_{A_{\ov{F}}}(L_{\ov{F}}) = Q_{\ov{F}}$. Since $r$ is
$\tau$-symmetric, we also have $r + \tau(r) = 2r \in A_{\ov{F}}^*$
(since $char(F) \neq 2$). Now setting $q = w^{-1} r$, we have $wq =
r$, and so $q$ satisfies the required hypotheses - i.e. $q \in
U(\ov{F})$.
\end{proof}

\begin{proof}[proof of theorem \ref{exp2}]
Let $L_1, L_2$ be subfields of $A$ of degree $m$, represented by
points $[L_1], [L_2] \in U(F) \subset \et_m(A)(F)$. By the previous
lemma, we may find a symplectic involutions $\sigma_1, \sigma_2$ such
that $L_i \subset Sym(A, \sigma_i)$. By \cite{BofInv}, proposition
2.7, there is an element $u \in Sym(A, \sigma)$ such that $\sigma_2 =
inn_u \circ \sigma_1$, where $inn_u$ denotes conjugation by
$u$. Define a morphism $\bA^1 \ra Sym(A, \sigma)$ by mapping $t$ to
$v_t = t u + (1 - t)$. Note that since $v_t \in Sym(A, \sigma)$, for
$t \in U$ we have that $\gamma_t = inn_{v_t} \circ \sigma_1$ is a
symplectic involution by \cite{BofInv}, proposition 2.7. Let
$Prp_{\sigma_i}$ be the Pfaffian characteristic polynomial on $Sym(A,
\sigma_i)$ (see \cite{BofInv}, page 19), which is a degree $m$
polynomial. Every element in $Sym(A, \sigma_i)$ satisfies the degree
$m$ polynomial $Prp_{\sigma_i}$, and further, there are dense open
sets of elements in $Sym(A, \sigma_i)$, $i = 1,2$ whose Pfaffians have
distinct roots (for example we may choose generators of $L_1$ and
$L_2$ respectively. Since $Sym(A, \sigma_i)$ is a rational variety and
$F$ is infinite, the $F$-points in $Sym(A, \sigma_i)$ are dense.  Note
also that since $[L_i] \in U(F)$, the restriction $[F(a)] \in U$ gives
a nonempty open condition on $Sym(A, \sigma_i)$.  Therefore, there is
an element $\alpha_1$ in $Sym(A, \sigma)$ such that the Pfaffians of
both $\alpha_1$ and $u \alpha_1 \in Sym(A, \sigma_2)$ have distinct
roots, and such that $[F(\alpha_i)] \subset U(F)$.

Let $\alpha_2 = u
\alpha_1$, and set $E_i = F(\alpha_i)$.

We will now show that the points $[L_1]$ and $[L_2]$ are R-equivalent
by first showing $[L_i]$ may be connected to $[E_i]$ by a rational
curve, and then showing $[E_1]$ and $[E_2]$ may also be connected by a
rational curve.

To connect $[E_1]$ and $[E_2]$, we define $\phi : \bA^1 \ra A$ via
$\phi(t) = v_t \alpha_1$. By construction, $E_1 = F(\phi(0))$ and $E_2
= F(\phi(1))$. Note also that $\phi(t) \in Sym(A, \gamma_t)$, and so
it satisfies the Pfaffian characteristic polynomial
$Prp_{\gamma_t}$. The condition that $Prp_{\gamma_t}(\phi(t))$ has
distinct roots gives an open condition on $t$ which is nontrivial
(e.g. $t = 0, 1$), and the condition that $[F(\phi(t))] \in U$ also
gives a nontrivial open condition, which together define a Zariski
dense open set of $\bA^1$. As in the proof of theorem
\ref{max_surject}, we may therefore obtain a rational morphism $\bA^1
\dra U \subset \et_m(A)$ via $t \mapsto [F(\phi(t))]$. It is easy to
check that this morphism sends $0$ to $[E_1]$ and $[1]$ to $[E_2]$.

Finally, to connect $[L_i]$ and $[E_i]$, choose a generator $\beta_i$
of $L_i$. Since both $\beta_i$ and $\alpha_i$ are
$\sigma_i$-symmetric, and since the $\sigma_i$-symmetric elements of
$A$ form a linear space, we may obtain a morphism $\bA^1 \ra Sym(A,
\sigma_i)$ by $t \mapsto t \beta_i + (1 - t) \alpha_i$. Since the
general element in the image of this morphism has distinct roots for
its Pfaffian, and generates a \'etale algebra in $U$, we obtain, as in
the proof of theorem \ref{max_surject} a rational morphism $\bA^1 \dra
\et_m(A)$ with $0 \mapsto [F(\alpha_i)] = [E_i]$ and $1 \mapsto
[F(\beta_i)] = [L_i]$.
\end{proof}

\section{Severi-Brauer flag varieties} \label{gsbv}

Suppose $A/F$ is a central simple algebra with $char(F)$ not dividing
$ind(A)$. Recall $V_d(A)^{(m)}_*$ denotes the open set in
$V_d(A)^{(m)}$ consisting of ideals which are linearly independent as
subspaces of $A$.

\begin{Lem} \label{ind_gsbv}
Suppose $A$ is a central simple algebra of index $i$ and $md =
i$. Then $ind(V_d(A)) = m$.
\end{Lem}
\begin{proof}
Without loss of generality, we may assume $F$ is prime to $p$ closed,
$p \neq char(F)$. By \cite{Blanchet}, $V_d(A)(L) \neq \emptyset$ if
and only if $ind(A_L) | d$. Therefore, it suffices to consider the
case that $A$ is a division algebra and that $d$ is a power of
$p$. Let $E \subset A$ be a maximal separable subfield. Since $F$ is
prime to $p$ closed and $p \neq char(F)$, $E$ has a Galois closure
which is a $p$-group and so it has subextensions of every size
divising $i = [E:F] = deg(A)$. In particular, $V_d(A)^{(m)}_*(F) =
\et_m(A)(F) \neq \emptyset$, and so $ind(V_d(A)) | m$.

On the other hand, suppose $V_d(A)(L) \neq \emptyset$ for some field
$L$. Since $ind(A_L) | d$, we may choose a maximal subfield $E \subset
A_L$ with $[E:L] | d$. Since $E$ splits $A$, it must contain a maximal
subfield of $A$, and so $deg(A) = i | [E:F] = [E:L] [L:F] |
d[L:F]$. Therefore, we have $m | [L:F]$, and in particular, $m |
ind(V_d(A))$, completing the proof.
\end{proof}

\begin{Lem} \label{open_index_gsbv}
Suppose $A$ is an $F$-central simple algebra of degree $n$ with index
$i$, and suppose $i = md$. If either $i$ is prime to $char(F)$ or $F$
is perfect, then $V_d(A)^{(m)}(F) = V_d(A)^{(m)}_*(F)$.
\end{Lem}
\begin{proof}
Suppose $x \in V_d(A)^{(m)}(F)$. Write $x \cong Spec(L)$ for $L/F$ a
degree $m$ \'etale extension. Choose a Galois extension $E/F$ with
group $G$ such that $L \otimes E \cong \oplus^m E$. Then $x$ gives a
collection of $m$ ideals $I_1, \ldots, I_m \subset A_E$ each of
reduced dimension $d$. Setting $I = \sum I_i$, note that the natural
$G$ action on $A_E$ restricts to an action on $I$, and hence by
descent, $I$ corresponds to an ideal $\ov{I} \subset A$ with
$rdim(\ov{I}) = rdim(I)$. In particular, since $ind(A) = i$, we have
$i | rdim(I)$. But this means the ideals $I_i$ are all linearly
independent and therefore $x$ corresponds to a point of
$V_d(A)^{(m)}_*(F)$ as claimed.
\end{proof}

\begin{Thm} \label{0_gsbv}
Let $X = V_{n_1, \ldots, n_k}(A)$. Let $$d = gcd\{n_1, \ldots, n_k,
ind(A)\}.$$ Then $X^{(ind(X))}$ is R-trivial if any of the following
conditions hold:
\begin{enumerate}
\item \label{d=1}
  $d = 1$,
\item \label{d=2}
  $d = 2$ and either $ind(A) | 4$ or $exp(A) | 2$,
\item \label{rel_prime}
  $d$ and $\frac{deg(A)}{gcd\{deg(A), d\}}$ are relatively prime.
\end{enumerate}
In particular, in each of these situations we have $A_0(X) = 0$.
\end{Thm}

\begin{proof}
Let $D$ be the underlying division algebra of $A$, and let $Y =
V_d(D)$. By theorem \ref{flag_reduction}, $Y$ and $X$ are stably
R-isomorphic. Therefore by proposition \ref{stably_rtriv}, it suffices
to show $Y^{(m)}$ is R-trivial. If we let $i = ind(A) = deg(D)$, then
we have $m = ind(Y) = i/d$ by lemma \ref{ind_gsbv}. By lemma
\ref{open_index_gsbv}, the inclusion $\et_m(D) \subset Y^{(m)}$ is
surjective on $F$-points. Therefore it suffices to show that
$\et_m(D)$ is R-trivial. But this follows from theorems \ref{maximal},
\ref{degree4} and \ref{exp2} respectively.
\end{proof}

\section{Involution varieties} \label{inv}

We assume in this section that the field $F$ has characteristic not
$2$.

\begin{Def}
Let $(A, \sigma)$ be an algebra with an involution (always assumed to
be of the first kind, either orthogonal or symplectic). We define the
radical of a right ideal $I \subset A$ to be $I \cap I^{\perp}$, where
$I^{\perp} = r.ann(\sigma(I))$. We say that a right ideal $I$ is
regular with respect to $\sigma$ if $rad(I) = 0$, or equivalently, $A
= I \oplus I^{\perp}$.
\end{Def}

We let $V_i(A)_{reg}$ be the subscheme of $V_i(A)$ consisting of
regular ideals. It is not hard to show that $V_k(A)_{reg}$ forms an
open subvariety of the generalized Severi-Brauer variety
$V_k(A)$. Hence, $(A, \sigma)$ has an regular ideal of reduced
dimension $k$ if and only if $ind(A) | k$.

We define the generalized involution variety $V_k(A, \sigma)$ to be
the subvariety of the Grassmannian representing the following functor
of points:

\begin{equation} \label{inv_def}
V_k(A, \sigma)(R) = 
\left\{I \in Gr(n^2 - nk, A)(R) \left| 
\begin{matrix}
I \text{ is a left ideal of } A_R \\
\text{and } \sigma(I) I = 0
\end{matrix}
\right.\right\}
\end{equation}

When $k = 1$, we write $V(A, \sigma)$ for $V_1(A, \sigma)$ and call
this the involution variety associated to $(A, \sigma)$.

\begin{Def}
Let $I$ be a right ideal of $(A, \sigma)$, and choose $l <
rdim(I)$. Define the subinvolution variety $V_l(I, \sigma)$ as the
variety representing the functor:
\begin{equation*}
V_l(I, \sigma)(R) = \{J \in V_l(A, \sigma)(R) | J \subset I\}
\end{equation*}
\end{Def}

The behavior of this variety depends on the ideal $I$ - in particular
on whether it is regular, isotropic or neither.

\begin{Thm} \label{subinv}
Suppose $(A, \sigma)$ is an algebra with involution. Let $I \subset A$
be an regular right ideal of reduced dimension $k$. Then there exists
a degree $k$ algebra with involution of the same type $(D, \tau)$
which is Brauer equivalent to $A$ and such that for any $l \leq k$, we
have:
\begin{equation*}
V_l(I, \sigma) = V_l(D, \tau)
\end{equation*}
\end{Thm}

\begin{proof}[proof of theorem \ref{subinv}]
By the fact that $I$ is regular, we may write $A = I \oplus I^\perp$,
and as in lemma \ref{idempotent}, $I = eA$ where $1 = e + f$, with $e
\in I, f \in I^{\perp}$. We set $D = eAe$. By \cite{Pie}, $D$ is
Brauer equivalent to $A$. By descent, one sees that $\sigma(e) = e$.
This implies that the involution $\sigma$ restricts to an involution
of $D$, and we denote this restriction by $\tau$.

To prove the theorem, we will construct mutually inverse maps (natural
transformations of functors) $\phi : V_l(I, \sigma) \ra V_l(D, \tau)$
and $\psi : V_l(D, \tau) \ra V_l(I, \sigma)$. For a commutative
Noetherian $F$-algebra $R$, and for $J \in V_l(I, \sigma)(R)$, we
define $\phi(J) = Je = eJe \subset D$. For $K \in V_l(D, \tau)$,
define $\phi(K) = KA$. It follows from an argument identical the one
in the proof of theorem \ref{sub_flag} that these are mutually inverse.
\end{proof} 

For an isotropic ideal, we have the following:

\begin{Lem} \label{isosubinv}
Suppose $(A, \sigma)$ is an algebra with involution. Let $I \subset A$
be an isotropic ideal of reduced dimension $k$. Then there exists a
degree $k$ algebra $D$ which is Brauer equivalent to $A$ such that for
any $l \leq k$,
\begin{equation*}
V_l(I, \sigma) = V_l(D)
\end{equation*}
\end{Lem}
\begin{proof}
This follows immediately from the fact that any ideal $J$ contained in
$I$ is automatically isotropic. Therefore, $V_l(I, \sigma) =
V_l(I)$. By theorem \ref{sub_flag}, we have $V_l(I) = V_l(D)$ as
claimed.
\end{proof}

\subsection{Orthogonal involution varieties}

\begin{Lem} \label{isoform_rtriv}
Suppose $V$ is a vector space space, and $q$ is an isotropic quadratic
form on $V$. Then the quadric hypersurface $C(q)$ is a rational
variety and any two $F$-points on $C(q) \subset \mathbb{P}(V)$ are
elementarily linked.
\end{Lem}
\begin{proof}
Since $q$ is isotropic, choose $p \in C(q)$. Consider the variety of
lines in $\mathbb{P}(V)$ passing through $p$. This is isomorphic to
$\mathbb{P}^{dim(V)-2}$, and hence is R-trivial. It is easy to see that the
rational morphism $\mathbb{P}^{dim(V)-2} \dra C(q)$ given by taking a line
through $p$ to its other intersection point with $C(q)$ is a
birational isomorphism, well defined off of the intersection of the
tangent space to $T_p C(q) \subset \mathbb{P}(V)$ with $C(q)$. In particular,
the $F$-points on $C(q)$ are infinite and dense. Now choose points
$p_1, p_2 \in C(q)(F)$. We may choose $p$ such that the rational
morphism defined above has both $p_i$ in its image by choosing $p$ in
the open complement of $T_{p_i} C(q)$. Using this map we may connect
$p_1$ and $p_2$ by a single rational curve by connecting their
preimages in $P^{dim(V)-2}$.
\end{proof}

\begin{Lem}
Suppose $(A, \sigma)$ is an algebra with orthogonal involution, and
let $X = V(A, \sigma)$. Then either $ind(X) = 1$ or $ind(X) =
max\{ind(A), 2\}$.
\end{Lem}
\begin{proof}
Consider the case where $ind(A) \leq 2$. In this case, we may choose
an ideal $I \subset V(A, \sigma)_{reg}$, and we have $V(I, \sigma)
\subset X$ is a subscheme which by descent is isomorphic to the
spectrum of a degree $2$ \'etale extension $E/F$. This means $ind(X)$
is $1$ or $2$. In particular, if $ind(A) = 2$, then $ind(X) = 2$ since
$V(I, \sigma) \subset X \subset V(A)$, and $ind(V(A)) = 2$, which
verifies the theorem in this case.

In the case $ind(A) > 2$, we may reduce to the case that $F$ is prime
to $2$ closed.  In particular, since $ind(V_2(A)) = ind(A)/2$ by lemma
\ref{ind_gsbv}, we may find a field extension $E/F$ of degree
$ind(A)/2$ such that $ind(A_E) = 2$ (note that $A_E$ is not split
since $[E:F] < ind(A)$). By the first case, $ind(X_E) = 2$, and so
there is a quadratic extension $L/E$ such that $X(L) \neq
\emptyset$. Therefore $ind(X) | ind(A)$. But since $X \subset V(A)$ and
$ind(V(A)) = ind(A)$, the reverse holds as well, and we have $ind(X) =
ind(A)$.
\end{proof}

\begin{Thm} \label{inv_rtriv}
Suppose $(A, \sigma)$ is a central simple $F$-algebra with orthogonal
involution and let $X = V(A, \sigma)$. If $F$ is prime to $2$-closed,
then $X^{(ind(X))}$ is R-trivial.
\end{Thm}
\begin{proof}
The case $ind(X) = 1$ follows immediately from lemma
\ref{isoform_rtriv}. 

If $ind(X) \geq 2$, we consider the morphism $f : X^{(2)}_* \ra
V_2(A)$ defined by taking a pair of ideals to their sum.  We note that
since a pair of $1$ dimensional subspaces are either equal or
independent, it follows by descent that $X^{(2)}_* = X^{(2)}$. Since
$ind(X) \neq 1$, every ideal $I \in V_2(A)(F)$ is either regular or
isotropic, since otherwise $rad(I)$ would be a point of
$X(F)$. Therefore it follows that the fiber over an ideal $I \in
V_2(A)(F)$ is $V(I, \sigma)^{(2)}$, which is either $Spec(E)^{(2)} =
Spec(F)$ for some quadratic \'etale $E/F$ (if $I$ is regular) or
$V(Q)^{(2)} = \et_2(Q)$ for some quaternion algebra $Q$ (if $I$ is
isotropic). In either case the fiber is nonempty (and R-trivial by
theorem \ref{maximal}).  Let $P = f^{-1}(V_2(A)_{reg})$. Since $f|_P$
is an isomorphism, we may regard $P$ as an open subvariety of
$V_2(A)$.

In the case $ind(X) = 2$, we have by corollary \ref{gsbv_rtriv} that
any two points in $P$ are elementarily linked. Therefore we may
conclude that $X^{(2)}$ is R-trivial if we can show that any point in
$X^{(2)}(F)$ is R-equivalent to one in $P(F)$. Let $\alpha \in
X^{(2)}(F)$ be arbitrary, let $J = f(\alpha)$, and choose a right
ideal $I \in V_2(A)_{reg}(F)$. By corollary \ref{gsbv_rtriv}, we may
find a morphism $\phi : \bP^1 \ra V_2(A)$ with $\phi(0) = J,
\phi(\infty) = I$. Since the generic point of $\bP^1$ maps into
$V_2(A)_{reg} \cong P \subset X^{(2)} \subset X^{[2]}$, we may lift
$\phi$ to a morphism $\psi : \bP^1 \ra X^{[2]}$ such that $\psi(0) \in
V(J, \sigma)^{[2]}$ and $f(\psi(\infty)) = I$. But since $ind(X) = 2$,
it follows from lemma \ref{index_hilb_open} that $V(J,
\sigma)^{[2]}(F) = V(I, \sigma)^{(2)}(F)$. In particular, since this
is an R-trivial variety, we find that $f^{-1}(I) = \psi(\infty) \sim_R
\psi(0) \sim_R \alpha$. Since $f^{-1}(I) \in P(F)$, $X^{(2)}$ is
R-trivial.

Suppose now that $i > 2$. By lemma \ref{connect_lem_new} it suffices
to show that $X^{(2)(i/2)}$ is R-trivial. Choose $\beta, \beta' \in
X^{(2)(i/2)}(F)$. In the case that $\beta, \beta' \in P^{(i/2)}(F)$,
the conclusion follows theorem \ref{exp2} since $P^{(i/2)} =
V_2(A)_{reg}^{(i/2)}$ and $(V_2(A)_{reg})^{(i/2)}_*(F) =
(V_2(A)_{reg})^{(i/2)}(F)$ and the fact that
$(V_2(A)_{reg})^{(i/2)}_*$ is an open subvariety of $\et_2(A)$.

Therefore we are done if we can show that for every $\beta \in
X^{(2)(i/2)}(F)$, there is a $\beta' \in P^{(i/2)}(F)$ with $\beta
\sim \beta'$. Given such a $\beta$, we may write $\beta =
\cH(\til{\beta})$ for some $\til{\beta} \in X^{(2)}(L)$, for $L/F$ a
degree $i/2$ field extension. By changing our focus to proving the
same thing for $\til{\beta}$, we may assume by lemma
\ref{new_transfer}, that $L = F$, $i = 2$, in which case we are done
by the argument in the $i = 2$ case.
\end{proof}

\begin{Thm} \label{0_orth}
Suppose $(A, \sigma)$ is a central simple algebra with orthogonal
involution. Then $A_0(V(A, \sigma)) = 0$.
\end{Thm}
\begin{proof}
This follows directly from theorems \ref{ch_thm} and \ref{inv_rtriv}.
\end{proof}

\subsection{Symplectic involution varieties}

Let $(A, \sigma)$ be an algebra with symplectic involution and index
at most $4$. Note that since every reduced dimension $1$ right ideal
is isotropic, the variety $V(A, \sigma)$ is the same as $V(A)$. We
therefore focus our attention on the first nontrivial case $V_2(A,
\sigma)$. 

\begin{Lem}
Let $X = V_2(A, \sigma)$ as above. Then $ind(X)$ is $1$ or $2$.
\end{Lem}
\begin{proof}
Without loss of generality, we may assume that $F$ is prime to
$2$-closed. Suppose $X(F) = \emptyset$. We must show that $X$ has a
point in a quadratic extension of $F$. Choose $I \in V_4(A)_{reg}$,
and consider the generalized subinvolution variety $V_2(I, \sigma)$
which is a closed subscheme of $X$. By theorem \ref{subinv}, $V_2(I,
\sigma) \cong V_2(D, \tau)$ for some degree $4$ algebra with
symplectic involution $\tau$. It therefore suffices to consider the
case that $deg(A) = 4$, and this follows from the following proposition.
\end{proof}

\begin{Prop} \label{symp_quadric}
Suppose $A$ is a degree $4$ algebra with symplectic involution
$\sigma$. Then $V_2(A, \sigma)$ is isomorphic to a quadric
hypersurface in $\bP^4$.
\end{Prop}
\begin{proof}
By \ref{degree4_pluker}, recall that the Pl\"uker embedding descends
to show $V_2(A)$ as a quadric hypersurface in $V(B)$ where $B$ is a
degree $6$ central simple algebra similar to $A^{\otimes2}$. In
particular, since $exp(A) | 2$, we have $V(B) \cong \bP^5$. 

At the separable closure, if we write $A = End(W)$, $B = End(\wedge^2
W)$, this corresponds to the Pl\"uker embedding $Gr(2, W) \hra
\bP(\wedge^2 W)$. The symplectic involution $\sigma$ is adjoint to a
form $\omega$ on $W$ which defines an element of $W^* \wedge W^* =
\cO_{\bP(\wedge^2)}(2)$, and the zeros off this element in $Gr(2, W)$
are exactly the isotropic subspaces. By descent, this corresponds to a
hyperplane in $\bP^5 = V(B)$, whose intersection with the embedded
$V_2(A)$ is $V_2(A, \sigma)$. Hence, by intersecting our quadric
$V_2(A)$ with an additional hyperplane, we obtain a quadric $V_2(A,
\sigma)$ in $\bP^4$ as claimed.
\end{proof}

\begin{Cor}
Suppose $X = V_2(A, \sigma)$ for an algebra $A$ of degree $4$. Then
$X^{(ind(X))}$ is R-trivial.
\end{Cor}
\begin{proof}
This follows from proposition \ref{symp_quadric} and theorem \ref{inv_rtriv}.
\end{proof}

\begin{Thm} \label{symp_rtriv}
Let $X = V_2(A, \sigma)$, and assume $F$ is prime to $2$-closed. Then
$X^{(ind(X))}$ is R-trivial.
\end{Thm}

\begin{proof}
We first consider the case $ind(X) = 1$. In this case, if we have two
points $I_1, I_2 \in X(F)$, note that the ideals $J \in X(F)$ such
that $J$ is linearly disjoint from the $I_i$'s form a dense open
subvariety $U \subset X(F)$. Since the group $Sp(A, \sigma)$ acts on
$X$ with dense orbits and is unirational, we may find an element $a
\in Sp(A, \sigma)(F)$ such that $a(I_1) \in U(F)$. In particular,
$U(F)$ is nonempty. Choose $J \in U(F)$, and let $V \subset X$ be the
dense open subscheme of right ideals of reduced dimension $2$ which
are linearly disjoint from $J$. We have a morphism $f : V \ra V_4(A)$
via $f(I) = J + I$. The image $Y$ of $f$ is open in the subvariety of
ideals $K \in V_4(A)$ with $J \subset K$. It follows from remark
\ref{rtriv_flag_fibers} that $Y$ is R-trivial. Note that if $K \in
Y(F)$, the fiber $f^{-1}(K)$ is open in $V_2(K, \sigma)$ which is
R-trivial and rational by lemma \ref{isoform_rtriv}. Therefore by
corollary \ref{unirat_fibers}, $V$ is R-trivial, which implies $I_1
\sim_R I_2$.

Now consider the case $ind(X) = 2$. Let $f : X^{(2)}_* \ra V_4(A)$ as
before. Given $\alpha \in X^{(2)}(F)$, I claim that $\alpha \sim_R
\alpha'$ for some $\alpha' \in X^{(2)}_*$ such that $f(\alpha')$ is a
regular ideal. To see this, we write $\alpha = \cH(\beta)$ for $\beta
\in X(L)$, $L/F$ a degree $2$ field extension. Since $Sp(A_L,
\sigma_L)$ is unirational (\cite{Borel}, thm 18.2) and acts on $X_L$
with dense orbits, we may choose $\psi : \bP^1_L \ra Sp(A, \sigma)$
such that $\psi(0) = id$, and $\alpha' = \psi(\infty)(\alpha)$ is in
the open set of elements such that $\alpha' \in X^{(2)}_*$ and
$f(\alpha')$ is regular. The path $\phi : \bP^1 \ra X(L)$ via $\phi(t)
= \psi(t)(\alpha)$ shows that $\alpha \sim_R \alpha'$ as claimed.

Let $P = f^{-1}(V_4(A)_{reg})$. We have reduced to showing that $P$ is
R-trivial. But by theorem \ref{subinv}, the fibers are of the form
$V_2(D, \tau)^{(2)}$ for $D$ an index $4$ algebra with symplectic
involution $\tau$. Since $V_2(D, \tau)$ is isomorphic to a quadric
hypersurface in $\bP^4$, by proposition \ref{symp_quadric}, $V_2(D,
\tau)^{(2)}$ is birational to the Grassmannian of projective lines in
$\bP^4$, and hence they are unirational of constant positive
dimension. Since they are also R-trivial by \ref{inv_rtriv}, we
conclude from corollary \ref{unirat_fibers} that $P$ is R-trivial as
desired.
\end{proof}

\begin{Thm} \label{0_symp}
Let $A$ be a central simple algebra with symplectic involution
$\sigma$ and and index at most $4$ and let $X = V_2(A, \sigma)$. Then
$A_0(X) = 0$.
\end{Thm}
\begin{proof}
This follows from theorems \ref{symp_rtriv} and \ref{ch_thm}.
\end{proof}

\bibliographystyle{alpha}
\bibliography{citations}

\def\cprime{$'$} \def\cprime{$'$} \def\cprime{$'$} \def\cprime{$'$}
\begin{thebibliography}{KMRT98}

\bibitem[Art82]{Ar:BS}
M.~Artin.
\newblock Brauer-{S}everi varieties.
\newblock In {\em Brauer groups in ring theory and algebraic geometry (Wilrijk,
  1981)}, pages 194--210. Springer, Berlin, 1982.

\bibitem[Bla91]{Blanchet}
Altha Blanchet.
\newblock Function fields of generalized {B}rauer-{S}everi varieties.
\newblock {\em Comm. Algebra}, 19(1):97--118, 1991.

\bibitem[Bor91]{Borel}
Armand Borel.
\newblock {\em Linear algebraic groups}, volume 126 of {\em Graduate Texts in
  Mathematics}.
\newblock Springer-Verlag, New York, second edition, 1991.

\bibitem[Bro]{Bro:BB}
Patrick Brosnan.
\newblock {On motivic decompositions arising from the method of
  Bialynicki-Birula}.
\newblock arXiv:math.AG/0407305.

\bibitem[CGM05]{CGM}
Vladimir Chernousov, Stefan Gille, and Alexander Merkurjev.
\newblock Motivic decomposition of isotropic projective homogeneous varieties.
\newblock {\em Duke Math. J.}, 126(1):137--159, 2005.

\bibitem[CT05]{CT}
J.-L. Colliot-Th{\'e}l{\`e}ne.
\newblock Un th\'eor\`eme de finitude pour le groupe de {C}how des
  z\'ero-cycles d'un groupe alg\'ebrique lin\'eaire sur un corps $p$-adique.
\newblock {\em Invent. Math.}, 159(3):589--606, 2005.

\bibitem[DI71]{DeIn}
Frank DeMeyer and Edward Ingraham.
\newblock {\em Separable algebras over commutative rings}.
\newblock Springer-Verlag, Berlin, 1971.

\bibitem[Ful98]{Ful:IT}
William Fulton.
\newblock {\em Intersection theory}, volume~2 of {\em Ergebnisse der Mathematik
  und ihrer Grenzgebiete. 3. Folge. A Series of Modern Surveys in Mathematics
  [Results in Mathematics and Related Areas. 3rd Series. A Series of Modern
  Surveys in Mathematics]}.
\newblock Springer-Verlag, Berlin, 1998.

\bibitem[Kar00]{Kar:Cell}
N.~A. Karpenko.
\newblock Cohomology of relative cellular spaces and of isotropic flag
  varieties.
\newblock {\em Algebra i Analiz}, 12(1):3--69, 2000.

\bibitem[KM90]{KarpMerk}
N.~A. Karpenko and A.~S. Merkur{\cprime}ev.
\newblock Chow groups of projective quadrics.
\newblock {\em Algebra i Analiz}, 2(3):218--235, 1990.

\bibitem[KMRT98]{BofInv}
Max-Albert Knus, Alexander Merkurjev, Markus Rost, and Jean-Pierre Tignol.
\newblock {\em The book of involutions}.
\newblock American Mathematical Society, Providence, RI, 1998.
\newblock With a preface in French by J.\ Tits.

\bibitem[Kol99]{Kol:RCLF}
J{\'a}nos Koll{\'a}r.
\newblock Rationally connected varieties over local fields.
\newblock {\em Ann. of Math. (2)}, 150(1):357--367, 1999.

\bibitem[KS03]{KolSza:RCFF}
J{\'a}nos Koll{\'a}r and Endre Szab{\'o}.
\newblock Rationally connected varieties over finite fields.
\newblock {\em Duke Math. J.}, 120(2):251--267, 2003.

\bibitem[KS04]{KraSa}
Daniel Krashen and David~J. Saltman.
\newblock Severi-{B}rauer varieties and symmetric powers.
\newblock In {\em Algebraic Transformation Groups and Algebraic Varieties},
  volume 132 of {\em Encyclopaedia Math. Sci.}, pages 59--70. Springer, Berlin,
  2004.

\bibitem[Mer95]{Merk:0cycles}
A.~S. Merkur{\cprime}ev.
\newblock Zero-dimensional cycles on some involutive varieties.
\newblock {\em Zap. Nauchn. Sem. S.-Peterburg. Otdel. Mat. Inst. Steklov.
  (POMI)}, 227(Voprosy Teor. Predstav. Algebr i Grupp. 4):93--105, 158, 1995.

\bibitem[Pie82]{Pie}
Richard~S. Pierce.
\newblock {\em Associative algebras}.
\newblock Springer-Verlag, New York, 1982.
\newblock Studies in the History of Modern Science, 9.

\bibitem[PSZ]{PSZ}
V.~Petrov, N.~Semenov, and K.~Zainoulline.
\newblock {Zero-cycles on a twisted Cayley plane}.
\newblock arXiv:math.AG/0508200.

\bibitem[Sam56]{Sam}
Pierre Samuel.
\newblock Rational equivalence of arbitrary cycles.
\newblock {\em Amer. J. Math.}, 78:383--400, 1956.

\bibitem[Swa89]{Swan}
Richard~G. Swan.
\newblock Zero cycles on quadric hypersurfaces.
\newblock {\em Proc. Amer. Math. Soc.}, 107(1):43--46, 1989.

\bibitem[SZ]{SeZa}
Nikita Semenov and Kirill Zainoulline.
\newblock {On motivic decompositions arising from the method of
  Bialynicki-Birula}.
\newblock preprint at http://www.mathematik.uni-bielefeld.de/LAG/man/189.html.

\bibitem[Tao94]{Tao}
David Tao.
\newblock A variety associated to an algebra with involution.
\newblock {\em J. Algebra}, 168(2):479--520, 1994.

\end{thebibliography}

\end{document}